\documentclass[11pt,a4paper]{article}
\usepackage[utf8]{inputenc}
\usepackage{amsmath}
\usepackage{amsfonts}
\usepackage{amssymb}
\usepackage{makeidx}
\usepackage{graphicx}
\usepackage{mathabx} %benötigt für schreibweise der rankine-hugoniot bedingung
\usepackage{dsfont} %benötigt für charakteristische funktion (1-funktion)
\usepackage[left=2cm,right=2cm,top=2cm,bottom=2cm]{geometry}
\usepackage{color}
\usepackage{verbatim} %füge kommentare ein durch \begin{comment} \end{comment}
\usepackage[nottoc]{tocbibind} %packt literaturverzeichnis in table of contents ohne table of contents selbst
\usepackage{hyperref} %um hyperlinks einzufügen
\usepackage{mathtools} % für inklusionspfeil
\usepackage{amsthm}
\usepackage{cancel}
\usepackage{todonotes}
\usepackage{float}
\usepackage{appendix}
\usepackage{emptypage} %removes headers and footers on empty pages

\numberwithin{equation}{section}

\newtheorem{lemma}{Lemma}[section]
\newtheorem{theorem}{Theorem}[section]
\newtheorem{definition}{Definition}[section]
\newtheorem{remark}{Remark}[section]
\newtheorem{proposition}{Proposition}[section]
\newtheorem{corollary}{Corollary}[section]
\begin{document}

\setlength{\parindent}{0pt} %damit nicht jeder neue absatz eingerückt ist

\begin{comment}

\section{Additional literature}

\end{comment}

\begin{center}
\Large\textbf{A VARIATIONAL APPROACH TO THE MODELING OF COMPRESSIBLE MAGNETOELASTIC MATERIALS} \\[10mm]
\large{BARBORA BENE\v{S}OVÁ}$^1$, \large{\v{S}ÁRKA NE\v{C}ASOVÁ}$^2$, \large{JAN SCHERZ}$^{1,2,3}$,\\ \large{ANJA SCHLÖMERKEMPER}$^3$ \\[10mm]
\end{center}

\begin{itemize}
\item[$^1$] Department of Mathematical Analysis, Faculty of Mathematics and Physics, Charles University in Prague, Sokolovská 83, Prague 8, 18675, Czech Republic
\item[$^2$] Mathematical Institute, Academy of Sciences, \v{Z}itná 25, Prague 1, 11567, Czech Republic
\item[$^3$] Institute of Mathematics, University of Würzburg, Emil-Fischer-Str. 40, 97074 Würzburg, Germany
\end{itemize}      

\bigskip

\begin{center}
\Large\textbf{Abstract} \\[4mm]
\end{center}

%We prove the existence of weak solutions to 
We analyze a model of the evolution of a (solid) magnetoelastic material. More specifically, the model we consider describes the evolution of a compressible magnetoelastic material with a non-convex energy and coupled to a gradient flow equation for the magnetization in the quasi-static setting. The viscous dissipation considered in this model induces an extended material derivative in the magnetic force balance. We prove existence of weak solutions %The proof of our existence result is 
based on De Giorgi's minimizing movements scheme, which allows us to deal with the non-convex energy as well as the non-convex state space for the deformation. In the application of this method we rely on the fact that the magnetic force balance in the model can be expressed in terms of the same energy and dissipation potentials as the equation of motion, allowing us to model the functional for the discrete minimization problem based on these potentials.

\section{Introduction}

In this article, we are devoted to the modeling of the evolution of (solid) \textit{magnetoelastic} materials. Such materials constitute a specific class of deformable ferromagnetic materials. As their mechanical and magnetic properties are interlinked, they may deform mechanically not only under the action of mechanic but also magnetic forces; this is known as the \textit{magnetostrictive} effect. This effect plays a crucial role in the construction of so-called magnetic actuators, cf.\ \cite{actuatorsandsensors,actuators}, which generate mechanical energy from changes in magnetic fields. Such actuators can for example be used as steering instruments in remote drug delivery (cf.\ \cite{yang}), where small capsule-shaped devices filled with medical drugs are navigated to the target area through the bloodstream of the human body. The \textit{inverse magnetostrictive} effect refers to the change of magnetic properties under mechanical stress and is also observed in nature. This effect, in turn, is exploited in sensors that measure mechanical stresses by converting them into alterations in their magnetic fields (cf.\ \cite{bienkowski2,actuatorsandsensors,grimes}), which find use in civil engineering: For instance, they help to prevent the collapse of buildings by detecting corrosion, fatigue or damage (see \cite{ausanio}). \\

In previous decades, continuum theories modelling electro-magneto-elastic behaviour have been extensively investigated and debated without having reached a consensus on, e.g., the general form of the stress tensor. We refer to \cite{emek, brown2,dorfmann, eringenmaugin, eriksen, kovetz, landaulifshitz,podio,schloemerkemper,tiersten,xiaobhattacharya} for comparisons, extensive discussions as well as further related literature. \\

In this article we focus on a quasi-static setting on the continuum level of a non-conducting material that is ferromagnetic and elastic potentially showing large deformations. (See below for related literature.) At the core of the energetic approach is a combination of the long-established and well-studied theories of \emph{hyperelasticity} (cf.\ \cite{kruzikroubicek,ogden}) and \textit{micromagnetics} (cf.\ \cite{brown,kruzikprohl}). In more detail, we assume that the state of the specimen is described by a deformation $\eta: \Omega_0 \to \mathbb{R}^3$ and a magnetization $M: \eta(\Omega_0) \to \mathbb{R}^3$. Here, $\Omega_0 \subset \mathbb{R}^3$ is the reference configuration, which is often chosen as the stress-free state of the specimen. Due to the fact that the magnetic and elastic properties are strongly interlinked in this materials, the overall magneto-elastic stored energy is of the following form (following e.g.\ \cite{brown,kruzikprohl,kruzikroubicek,elasticity})
\begin{equation}
\int_{\Omega_0} \phi(\nabla \eta, \nabla^2 \eta) + \psi(\nabla \eta, M \circ \eta) \, dX + \int_{\eta(\Omega_0)} A |\nabla_x M|^2-\frac{\mu}{2} M\cdot (H(M, \eta)) + \varphi(M) \, dx.
\label{energ-intro}
\end{equation}
Here, $\phi$ is the purely elastic part of the energy, $\psi$ models the magnetostriction, $A |\nabla_x M|^2$ is the exchange energy (notice that the subscript "$x$" is added to indicate that the gradient in the actual configuration is meant), $H$ is the stray field. Finally, $\varphi$ is a purely magnetic part of the energy, which can, e.g., model a saturation constraint on the magnetization. We will comment on the individual parts of the energy more thoroughly in Section \ref{modeldescriptionmagnetoelasticity}. However, we wish to stress that due to physical reasons and modeling consistency $\phi$ cannot be a convex function. Indeed, as the considered material cannot be infinitely compressed, the energy needs to blow up as the determinant of the deformation gradient tends to zero. This is expressed in the term \eqref{2935} in the elastic part of the energy below and causes our energy potential to be non-convex.\\

\begin{figure}
\centering
\includegraphics[scale=0.65]{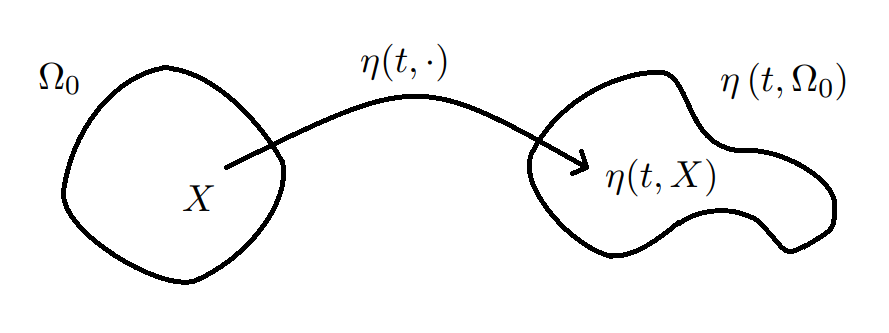} %[width=16cm,height=11cm]{solidinfluidnew}
%\hspace*{0.3cm}
\caption{A deformation $\eta(t,\cdot)$ mapping a (magneto)elastic material from its reference configuration $\Omega_0$ to its current configuration $\eta (t,\Omega_0)$.} \label{figure1}
\end{figure}

The modeling of magnetoelasticity can be traced back to the middle of the 20th century, cf., e.g., \cite{brown2,tiersten}. With regards to the mathematical analysis of magnetoelasticity, the most fundamental approach lies in the study of the static case, with the study of the existence of minimizers to the stored energy as given above: Efforts in this direction were made and the existence of such minimizers was proved in both the two-dimensional setting, see e.g.\ \cite{desimonedolzmann}, and the three-dimensional setting, see e.g.\ \cite{desimonejames,jameskinderlehrer,kruzikstefanellizeman,rybkaluskin}. In the latter work, the analysis was further extended to the quasistatic case, i.e.\ the evolutionary setting with inertial effects neglected. This was achieved using a time-incremental minimization scheme, also known as De Giorgi's minimizing movements method, cf.\ \cite{degiorgi}: The problem was discretized in time, and the existence of a solution at each discrete time could be concluded from the investigation of the steady-state setting. In that work, the evolution of the magnetization is assumed to be dissipative and rate-independent. Generalizations to the compressible case were achieved later in \cite{barchiesi} in the static setting and in \cite{bresciani} in both the static and the quasistatic setting. A further extension of the existence result in \cite{barchiesi} to the setting of deformations in the class of Orlicz-Sobolev spaces is available in \cite{stroffolini}.

In order to also take into account inertial effects in the mathematical analysis, some works restrict themselves to the setting of small strains; see e.g.\ \cite{carbou,chipot,ellahiani}.

In recent years quite some progress was made for models entirely phrased in Eulerian coordinates, see e.g.\ \cite{johannesarticle,francescojoshua,johannesthesis,joshuapaper,martinanja,jiang,joshuathesis,zabensky,zhao1,zhao2}. In that case, the modeling is more similar to fluids and the state of the specimen is described by its velocity, the magnetization, and the deformation tensor, i.e.\ a matrix that represents the deformation gradient in the actual configuration. Even if that approach allows one to consider an evolution that includes inertial effects as well as the full Landau-Lifschitz-Gilbert equation for the evolution of the magnetization, the drawback is that the shape of the sample is usually fixed throughout the evolution and the deformation tensor is heavily regularized so that it can no longer be identified with the deformation gradient in the reference configuration (although the latter can be avoided by regularizing the dissipative stress as in \cite{roubiceknew}).

Furthermore, magnetoelasticity has also been studied in combination with thermodynamics. The existence of weak solutions has been proven for models of thermally and electrically conducting magnetoelastic materials in \cite{roubicektomassetti} and thermally conducting magnetoelastic materials under diffusion in \cite{tr3}. In \cite{tr1}, the existence of weak solutions was achieved for a model of a thermally conducting magnetoelastic material - undergoing a ferro-to-paramagnetic phase transition - in Eulerian coordinates.\\

In the present work, we take a Lagrangian point of view and we consider a viscous evolution for the magnetization which utilizes a generalized time derivative that models the transport of the magnetization via the deformation as in \cite{johannesarticle,johannesthesis} instead of a rate-independent evolution. Thus, we work on the level of force balances, i.e.\ the level of partial differential equations. Nevertheless, our approach to proving the existence of weak solutions to the obtained partial differential equations will still be variational. Indeed, we propose a variant of the minimizing movements (or time-incremental) scheme fitted to the viscous magnetoelastic contribution that we consider. 

The minimizing movements scheme provides a decisive advantage over other approaches: It allows us to deal with the constraint that only injective functions can be accepted as deformations as well as the non-convex energy functional. Indeed, fixed point arguments, which are often used for solving coupled systems of partial differential equations, usually rely on convexity of the state space and thus do not constitute an option in our setting. This in particular excludes the classical Galerkin method as a possible approach. But also the Rothe method, in which the system is decoupled by a time discretization and the discrete equations are solved directly instead of via minimization is out of question. This is because such an approach requires convexity of the energy for a discrete form of the chain rule, which in turn is necessary for the deduction of a priori bounds of the discrete solution (see also Remark \ref{minimizingmovementsandnonconvexity}). In the minimizing movements scheme, however, such bounds are obtained naturally by comparing, at each discrete time, the value of the minimized functional in its minimizer against the same functional in the minimizer from the previous discrete time.

In order to be able to apply this method, however, we must scrutinize a suitable dissipation potential from which the gradient flow in \eqref{2938} below can be deduced.

Our method bears further the advantage, that it is closely fitted to the implementation in \cite{variationalapproachtofsi}. Therefore, following the approach in \cite{variationalapproachtofsi}, an extension of our results to the case involving inertia seems perfectly feasible via the introduction of a two time-scale (cf.\ Remark \ref{inertia-ext}).

The article is organized as follows: We present the physical model and introduce some additional mathematical notation in Section \ref{modelandnotation}. In Section \ref{approximatesystemmagnetoelastic}, we outline our approach to the proof of the existence of weak solutions to this model, which constitutes our main result. In Section \ref{weaksolmagnetoel}, we explicitly state our definition of weak solutions and our main result and finally, in Section \ref{proofsection}, we give the full proof of the result.

%The article is organized as follows: We present the model in Section \ref{modeldescriptionmagnetoelasticity} and consider the corresponding weak formulation and our main result in Section \ref{weaksolmagnetoel}. In Section \ref{approximatesystemmagnetoelastic} we introduce our set up for the minimizing movements scheme for the proof of this result. The proof is carried out subsequently in Section \ref{proofsection}.

We also remark that the same model as here has been analyzed in Chapter 5 of the PhD thesis \cite{thesis} of the third author, of which the present article constitutes a self-contained version.

\section{Model and notation} \label{modelandnotation}

\subsection{Model}
\label{modeldescriptionmagnetoelasticity}

In this section we present the model under consideration in this paper. We assume the reference configuration of the material to be a bounded domain $\Omega_0 \subset \mathbb{R}^3$, the boundary of which is divided into a free part $N$ and a part $P$ with a prescribed Dirichlet boundary condition,
\begin{align}
N \subset \partial \Omega_0,\quad \quad P:= \partial \Omega_0\setminus N. \nonumber
\end{align}

The deformation of the material is described by a mapping $$\eta : (0,\infty)\times \Omega_0 \rightarrow \mathbb{R}^3,$$ cf.\ Figure~\ref{figure1}. The current/actual configuration $\Omega_0 (t)$ of the material at any time $t \in (0,\infty)$ is expressed as
\begin{align}
\Omega (t) := \eta \left(t, \Omega_0 \right) \nonumber
\end{align}
and we define the time-space domain
\begin{align}
Q := \left\lbrace (t,x) \in (0,\infty)\times \mathbb{R}^3:\ x \in \Omega (t) \right\rbrace. \nonumber
\end{align}

In addition to the deformation, the state of the sample is characterized by its magnetization $M:Q \rightarrow \mathbb{R}^3$. The magnetization $M$ is most naturally defined in the actual configuration. Nonetheless, we may assign to $M$ its counterpart in Lagrangian coordinates denoted $\tilde{M}: (0,\infty)\times \Omega_0 \rightarrow \mathbb{R}^3$. We assume that, as the material undergoes a deformation, the modulus of the magnetization vector changes accordingly, which is mathematically expressed trough the relation
\begin{align}
\tilde{M} \left( t, X \right) = \det \left( \nabla \eta(t,X) \right) M (t,\eta(t,X)), \quad \quad \forall t \in (0,\infty) \text{ and } X \in \Omega_0. \label{magentization-trf}
\end{align}
We point out that this is a purely kinematic modeling assumption which is based on the corresponding assumption in the static case in \cite{desimonedolzmann, kinderlehrer}. We assume that this condition generalizes in an unchanged manner also to the evolutionary setting and holds for all times $t$ (as already done in e.g. \cite{johannesarticle,johannesthesis}). We adopt this approach as it seems to represent the simplest possible choice of generalization to the dynamic setting.

%Following \cite{desimonedolzmann, kinderlehrer}, we assume 
%\begin{align}
%\tilde{M} \left( t, X \right) = \det \left( \nabla \eta(t,X) \right) M (t,\eta(t,X)), \quad \quad \forall t \in (0,\infty) \text{ and } X \in \Omega_0, \label{magentization-trf}
%\end{align}
%which models that under the deformation of the material the modulus of the magnetization vector changes accordingly. The gradient $\nabla$ in \eqref{magentization-trf} is taken with respect to the Lagrangian variable $X$.

\begin{remark}[Transport of the magnetization]
By assuming formula \eqref{magentization-trf}, we assume that by elastic deformation the magnetization may be stretched but not rotated.  However, more complicated transformation rules for the magnetization vector could be assumed. For example, one could require $\tilde{M} \left( t, X \right) = \det \left( \nabla \eta(t,X) \right) (\nabla \eta)^{-1} M (t,\eta(t,X))$, which would conserve the flux of the magnetization through closed surfaces, as in, e.g.\ \cite{antman,bresciani,rogers,tr3}. In such a case, the form of the material derivative in \eqref{2938} also needs to be changed accordingly.
\end{remark}

The evolution of the state variables is given by the following abstract balance equations
\begin{align}
-\operatorname{div} \sigma (\eta, \tilde{M}) &=\varrho f \circ \eta \quad \quad \text{in } (0,\infty) \times \Omega_0, \label{2949} \\
\partial_t M + \left( v \cdot \nabla_x \right) M + \left( \nabla_x \cdot v \right)M &= -H_\text{eff} \quad \quad \text{in } Q, \label{2948}
\end{align}
supplemented by the boundary condition
\begin{align}
\eta (t) = \gamma \quad \quad \text{on } P \label{2947}
\end{align}

and suitable initial conditions. Here, \eqref{2949} is the balance of momentum with the inertial forces already neglected and $\sigma$ is the first Piola-Kirchhoff stress tensor of the material. Assuming a Kelvin-Voigt rheology, it can be determined via
\begin{align}
\operatorname{div} \sigma = D_1 \tilde{E} \left(\eta, \tilde{M} \right) + D_2 \tilde{R} \left( \eta, \partial_t \eta, \partial_t \tilde{M} \right), \label{2933}
\end{align}
where $D_1$ and $D_2$ denote the Fréchet derivatives of the energy $\tilde{E}$ (see \eqref{2999}) and the dissipation $\tilde{R}$ (see \eqref{2936}) with respect to the first and second argument, respectively. Finally, $\rho$ is the density of the material and $f$ is a body force applied in the actual configuration.

Equation \eqref{2948} is a gradient flow assuming the extended time derivative $D_t M$ 
\begin{align}
D_t M := \partial_t M + \left( v \cdot \nabla_x \right) M + \left( \nabla_x \cdot v \right) M, \quad \text{with} \quad v \left(t, x \right) = \left[\partial_t \eta (t)\right] \left( \eta^{-1}(t,x) \right),  \quad (t,x) \in Q. \label{2938}
\end{align}
The right-hand side of the equation \eqref{2948} is given by the effective magnetic field determined as an external magnetic field $H_{\text{ext}}$ and a Fréchet derivative of the total energy in the actual configuration (see Remark \ref{energy-Eulerian} for a precise form); i.e.
$$
H_\text{eff}:=  -\mu H_\mathrm{ext} + D_2 E \left(\eta, M \right).
$$

\begin{remark}
Let us note that the form of the extended derivative in \eqref{2938} is given by \eqref{magentization-trf}. Indeed, differentiating that relation with respect to time yields 
\begin{align*}
\partial_t \tilde M = \mathrm{det}(\nabla \eta) \bigg((\nabla \eta)^{-\top}:(\partial_t \nabla) \eta M +  \big( \partial_t M + (\partial_t \eta \cdot \nabla_x) M \bigg),
\end{align*}
where the expression in the bracket corresponds to $D_t M$ in Eulerian coordinates. Thus, the extended time derivative follows the same transformation rule as the magnetization itself. 
\end{remark}

We proceed by recalling the underlying energy from \eqref{energ-intro}, making it more precise, commenting on the individual terms, and rewriting it in its purely Lagrangian and Eulerian forms. We start with the Lagrangian one and set 
\begin{align}
&\tilde{E} \left( \eta, \tilde{M} \right) := \left\{
\begin{array}{ll} \tilde{E}_{\text{el}}\left(\eta \right) + \tilde{E}_{\text{mag}}\left( \eta, \tilde{M} \right)  \quad \quad &\text{if } \det \left(\nabla \eta \right) > 0 \ \text{a.e. in } \Omega_0, \\ + \infty \quad \quad &\text{otherwise}, \end{array} \right. \label{2999}
\end{align}
where
\begin{align}
\tilde{E}_{\text{el}}\left(\eta \right) :=& \int_{\Omega_0} W \left( \nabla \eta \right) + \frac{1}{\left( \det \left( \nabla \eta \right) \right)^a} + \frac{1}{q} \left| \nabla^2 \eta \right|^q \ dX, \label{4039} \\
\tilde{E}_{\text{mag}}\left( \eta, \tilde{M} \right) :=& \int_{\Omega_0} \tilde{\Psi}\left( \nabla \eta, \tilde{M} \right) - \frac{\mu}{2} \tilde{M} \cdot \left(H \left[ \tilde{M}, \eta \right](\eta) \right)\nonumber + A \left| \nabla \left( \frac{1}{\det \left( \nabla \eta \right)} \tilde{M} \right) \left( \nabla \eta \right)^{-1} \right|^2 \det\left( \nabla \eta \right) \\ &+ \frac{1}{4\beta^2} \left( \left| \frac{1}{\det \left( \nabla \eta \right)} \tilde{M} \right|^2 - 1 \right)^2 \det \left( \nabla \eta \right)\ dX. \label{4041}
\end{align}

First, let us notice that we set $\tilde{E}(\eta, \tilde{M}) = + \infty$ if the determinant of the deformation gradient takes non-positive values; this assures that the deformation is orientation preserving. In the elastic energy $\tilde{E}_{\operatorname{el}}$ the quantity $W: \mathbb{R}^{3 \times 3} \rightarrow \mathbb{R}_0^+$ denotes the elastic energy density. To the elastic energy density, we add the expression
\begin{align}
\frac{1}{\left( \det \left( \nabla \eta \right) \right)^a} + \frac{1}{q} \left| \nabla^2 \eta \right|^q, \label{2935}
\end{align}
in which we assume that $q > 3$ and $a > \frac{3q}{q-3}$. 

Here, the first term has both a physical and a mathematical meaning. It causes the elastic energy to blow up if the determinant of the deformation gradient becomes small; i.e.\ it models the resistance of the material to infinite compression. We also note that this notably makes the elastic part of the energy non-convex. Together with the second term and the assumed growth, we will even obtain a strict lower bound of this determinant (see Lemma \ref{boundawayfromzero} in the appendix following \cite{healeykromer}). The mathematical significance of the latter term consists in guaranteeing the $C^1$-regularity of the deformation. This in combination with the bound of the determinant of the deformation gradient away from zero allows for unproblematic transformations between the reference configuration and the current configuration, which is needed for the model to be well-posed due to the Eulerian-Lagrangian viewpoint taken in the modeling.

 In the micromagnetic energy $\tilde{E}_{\operatorname{mag}}$ the quantity $\tilde{\Psi}: \mathbb{R}^{3\times 3}\times \mathbb{R}^3 \rightarrow \mathbb{R}_0^+$ represents the anisotropy energy density (cf.\ \cite[Section 2.1]{desimone}) that models the magnetostriction. The stray field $H[\tilde{M},\eta]: (0,\infty)\times \mathbb{R}^3 \rightarrow \mathbb{R}^3$ is defined in the current configuration through the reduced set of Maxwell's equations as follows (see e.g.\ \cite[Section 1]{gar07}): 
 $$
 H[\tilde{M},\eta] := - \nabla \phi, 
 $$ with 
\begin{align}
\Delta_x \phi (t,\cdot) = \operatorname{div} M(t,\cdot) \quad \text{in } \mathbb{R}^3,\quad \quad \left\{
\begin{array}{ll} \left( \phi^{\operatorname{int}}(t,\cdot) - \phi^{\operatorname{ext}}(t,\cdot) \right) = 0 \quad &\text{on } \partial \Omega (t), \\ \left(\nabla_x \phi^{\operatorname{int}}(t,\cdot) - \nabla_x \phi^{\operatorname{ext}}(t,\cdot) \right)\cdot \operatorname{n} = - M \cdot \operatorname{n} \quad &\text{on } \partial \Omega (t), \end{array}
\right. \label{2932}
\end{align}
where $M(t,\cdot)$ has been extended by $0$ outside of $\Omega (t)$, $\phi^{\operatorname{int}}(t,\cdot)$ and $\phi^{\operatorname{ext}}(t,\cdot)$ denote the traces of $\phi (t,\cdot)$ to the interior and the exterior of $\Omega (t)$ respectively and $\operatorname{n}$ denotes the outer unit normal vector on $\partial \Omega (t)$. The constant $\mu > 0$ denotes the magnetic permeability. The expression
\begin{align}
A \left| \nabla \left( \frac{1}{\det \left( \nabla \eta \right)} \tilde{M} \right) \left( \nabla \eta \right)^{-1} \right|^2 \det\left( \nabla \eta \right), \nonumber
\end{align}
in which $A>0$ is called the exchange stiffness constant, represents the exchange energy. The magnetization of a ferromagnet has the tendency to align in a constant direction; this expression penalizes the deviation from such a behavior, cf.\ \cite[Section 3.2.2]{hs98}, \cite[§39]{landaulifshitz}. Another penalization term in the micromagnetic energy is given by
\begin{align}
\frac{1}{4\beta^2} \left( \left| \frac{1}{\det \left( \nabla \eta \right)} \tilde{M} \right|^2 - 1 \right)^2 \det \left( \nabla \eta \right), \label{2934}
\end{align}

cf.\ for example \cite{chipot}. For a small value of the constant $\beta>0$ this term penalizes the magnitude of the magnetization (in the current configuration) taking values away from one. In general, the magnetization of a ferromagnetic material is considered to be constant in magnitude for constant temperature and equal to the saturation magnetization (cf.\ \cite{desimone}, \cite[Section 9.2]{eringenmaugin}), which for simplicity can be set equal to one using a scaling argument. We point out that the quantity \eqref{2934} further implies non-convexity of the energy $\tilde{E}$ also with respect to the magnetization. 

\begin{remark}
\label{energy-Eulerian}
By a change of variables and by using relation \eqref{magentization-trf}, we may rewrite the energy also into the actual configuration as      
\begin{align}
&E \left( \eta, M \right) := \left\{
\begin{array}{ll} E_{\operatorname{el}}(\eta) + E_{\text{mag}} \left(\eta, M \right) = \tilde{E}\left(\eta, \tilde{M} \right)  \quad \quad &\text{if } \det \left(\nabla \eta \right) > 0 \ \text{a.e. in } \Omega_0, \\ + \infty \quad \quad &\text{otherwise} \end{array} \right. \nonumber
\end{align}

with
\begin{align}
E_{\operatorname{el}}(\eta):=& \int_{\eta(\Omega_0)} \frac{1}{\operatorname{det} \left( \left[\nabla \eta \right] \left(\eta^{-1} \right) \right)} W \left( \left[\nabla \eta \right] \left(\eta^{-1} \right) \right) + \frac{1}{\left(\operatorname{det} \left( \left[\nabla \eta \right] \left(\eta^{-1} \right) \right) \right)^{a+1}} \nonumber \\
&+ \frac{1}{q\operatorname{det} \left( \left[\nabla \eta \right] \left(\eta^{-1} \right) \right)} \left| \left[ \nabla^2 \eta \right]\left(\eta^{-1} \right) \right|^q\ dx, \nonumber \\
E_{\text{mag}} \left(\eta, M \right) :=& \int_{\eta (\Omega_0)} \frac{1}{\operatorname{det} \left( \left[\nabla \eta \right] \left(\eta^{-1} \right) \right)}\tilde{\Psi}\left( \left[\nabla \eta \right] \left(\eta^{-1} \right), \operatorname{det} \left( \left[\nabla \eta \right] \left(\eta^{-1} \right) \right) M \right) \nonumber \\
& - \frac{\mu}{2} M \cdot H \left[ \det \left(\nabla \eta \right) M, \eta \right] + A \left| \nabla M \right|^2 + \frac{1}{4\beta^2} \left( \left| M \right|^2 - 1 \right)^2\ dx. \nonumber
\end{align}
\end{remark}

As the system is dissipative, we also have an underlying dissipation which is best defined in Eulerian coordinates: 
\begin{align}
R \left( \eta, v, D_t M \right) := \int_{\eta (\Omega_0)} \nu \left| \nabla_x v \right|^2 + \frac{1}{2} \left| D_t M \right|^2\ dx = \tilde{R} \left(\eta, \partial_t \eta, \partial_t \tilde{M} \right). \label{alternativedissipation}
\end{align}
Here, the first term models the viscous dissipation in the solid via a Newtonian-fluid like dissipation potential. The second term is just the square of the generalized time derivative given in \eqref{2938}, which reflects a general principle that a gradient flow evolution corresponds to a viscous one with a quadratic dissipation potential. This allows us to rewrite \eqref{2948} as
\begin{align}
D_2E \left(\eta, M \right) + D_3R \left(\eta, v, D_tM \right) - \mu H_{\operatorname{ext}} = 0 \quad \quad \text{in } Q, \label{alternativemagforcebalance}
\end{align}
which shows the underlying variational structure of the model that we will exploit in our proof of the existence of its weak solutions.

\begin{remark}
Similarly, as in the case of the energy, the dissipation potential can be also expressed in Lagrangian coordinates via 
\begin{align}
\tilde{R} \left( \eta, \partial_t \eta, \partial_t \tilde{M} \right) := \int_{\Omega_0} \nu \left| \nabla \partial_t \eta \left(\nabla \eta \right)^{-1} \right|^2 \det \left( \nabla \eta \right) + \frac{1}{2} \left| \frac{1}{\det \left(\nabla \eta \right)} \partial_t \tilde{M} \right|^2 \det \left( \nabla \eta \right)\ dX. \label{2936}
\end{align}    
\end{remark}

\begin{remark}[Extension to inertial problems]
\label{inertia-ext}
Having an underlying variational structure, it is feasible to expect that extending the equations to the case where inertia is not neglected is well possible. In fact, in \cite{variationalapproachtofsi} a general method for extending the minimizing movement approach from quasi-static to inertial problem has been proposed.
Indeed, the idea is to discretize $\partial_t^2 \eta$ via a difference quotient and to solve what was called a \emph{time-delayed problem} in \cite{variationalapproachtofsi} that reads
\begin{align} \label{eq:introTimeDelayed}
\mathrm{div}\left(\sigma(\eta(t), \tilde M(t))\right) + \varrho f\circ \eta(t)  = \frac{\partial_t \eta(t) - \partial_t \eta(t-h)}{h}
\end{align}
for some fixed $h$. After this discretization of the second time derivative, the problem becomes mathematically essentially equivalent to ours, so it only remains to pass to the limit $h \to 0$.
\end{remark}

\begin{remark}[Landau-Lifschitz-Gilbert equation]
Modelling the evolution of the magnetization via a gradient flow as in \eqref{2948} is not uncommon (see \cite{johannesthesis,sourav,zabensky,zhao1,zhao2}), but one might instead consider the LLG-equation written in our context as follows:
$$
D_t M = -\alpha M \times H_\mathrm{eff} - \beta M \times M \times H_\mathrm{eff}, 
$$
with $\alpha, \beta$ constants. 
Including this equation makes the analysis more difficult, but in a purely Eulerian setting has been considered e.g.\ in \cite{johannesarticle,johannesthesis}, wherein the existence of weak solutions was proved - subject to a smallness condition on the initial data - in the case of a strongly convex elastic energy density, a heavily regularized deformation tensor and under neglect of the stray field. In \cite{joshuapaper,joshuathesis}, this result was generalized to the setting of an only convex elastic energy density and with the stray field included and supplemented by the local-in-time existence of strong solutions and weak-strong uniqueness. In two-dimensional periodic domains, the global-in-time existence of weak solutions could even be proved without any smallness condition, cf.\ \cite{francescojoshua}. Global-in-time existence without a smallness assumption could also be proved for dissipative solutions in both 2D and 3D (cf.\ \cite{martinanja}). For strong solutions in 3D, local-in-time existence and uniqueness was achieved in \cite{du}. Moreover, local-in-time existence as well as (under a smallness assumption on the initial data) global-in-time existence of classical solutions in $\mathbb{R}^2$ and $\mathbb{R}^3$ in the case of a non-regularized deformation tensor and the elastic energy density $W(\theta) = \frac{|\theta|^2}{2}$ was shown in \cite{jiang}.
\end{remark}

\begin{remark}[Ferromagnetic and martensitc materials]
As an example for a material exhibiting a strong magnetostrictive effect, we mention the alloy \emph{NiMnGa}. This alloy falls into the class of so-called shape-memory alloys that exhibit an extra-ordinary stress-strain response,  see, e.g., \cite{bhattacharya} for an introduction. This is due to a diffusionless phase transformation going on in the material that changes the symmetry of the crystalline lattice; in particular, the material can then form a microstructure of its variants in the lower symmetric \emph{martensite} phase. The case of \emph{NiMnGa} is then special in the sense that it is not only martensitic but also ferromagnetic;  specific magnetizations are linked to particular variants so that a change of microstructure can be triggered by an applied magnetic field \cite{desimonejames}. This effect can, for example, be used to construct micropumps (see \cite{bhattacharyajames,micropumps}).
\end{remark}

\subsection{Notation} \label{notationsection}

Before we turn to the mathematical analysis of the model introduced in Section \ref{modeldescriptionmagnetoelasticity}, we introduce some additional mathematical notation. Since we deal with a domain which is entirely occupied by a magnetoelastic material, we are in particular in need of some notation related to the deformability of such a domain. 

In order to avoid self-penetration of the magnetoelastic body, we restrict ourselves to deformations injective on the interior of $\Omega_0 \subset \mathbb{R}^3$, or more precisely, deformations lying in the non-convex set $\mathcal{E}$ given by
\begin{align}
\mathcal{E} := \left\lbrace \eta \in W^{2,q}\left(\Omega_0; \right):\ \tilde{E}_{\operatorname{el}}\left( \eta \right) < \infty,\ \left| \eta \left(\Omega_0 \right) \right| = \int_{\Omega_0} \det \left(\nabla \eta \right) \ dX,\ \eta|_P = \gamma \right\rbrace, \label{statespace}
\end{align}

wherein $\tilde{E}_{\operatorname{el}}$ denotes the elastic energy defined in \eqref{4039} and $\gamma : P \rightarrow \mathbb{R}^3$ is an injective boundary deformation such that it can be extended to $\Omega_0$ as a deformation with finite elastic energy. The identity
\begin{align}
\left| \eta \left(\Omega_0 \right) \right| = \int_{\Omega_0} \det \left(\nabla \eta \right) \ dX \label{4016}
\end{align}

in the definition of the set $\mathcal{E}$ is called the Ciarlet-Nečas condition. This condition, which has been introduced in \cite{ciarletnecas}, can be interpreted as follows: Provided that $\Omega_0$ is of class $C^{0,1}$, any local homeomorphisms $\eta \in C^1(\overline{\Omega_0})$ satisfying this condition is in fact a global homeomorphism in $\Omega_0$, i.e.\ $\eta$ is injective except for possibly on the boundary $\partial \Omega_0$. In particular, this holds true for any $\eta \in \mathcal{E}$: Indeed, from the Morrey embedding we know that $\mathcal{E} \subset C^1(\overline{\Omega_0})$. Further, as $\eta \in \mathcal{E}$ satisfies $\tilde{E}_{\text{el}}(\eta) < \infty$, Lemma~\ref{boundawayfromzero} in the appendix implies that $\det (\nabla \eta) > 0$ in $\overline{\Omega_0}$. Consequently, by the inverse function theorem, $\eta$ constitutes a local and hence a global homeomorphism. 

The Ciarlet-Nečas condition can further be expressed as a variational formulation of friction-less contact. Indeed, this was shown for smooth deformations in the static case in \cite{ciarletnecas} and later generalized to our setting in \cite{antonin1,kromerroubicek,frictionlesscontact1}.

\begin{remark} \label{interiorpoints}
The set $\mathcal{E}$ constitutes a closed subset of the affine function space
\begin{align}
\left\lbrace \eta \in W^{2,q}\left(\Omega_0 \right):\ \eta|_P = \gamma \right\rbrace. \nonumber
\end{align}

The interior and the boundary of $\mathcal{E}$ can be characterized in the following way: For any $\eta \in \mathcal{E}$ it holds that
\begin{align}
\eta \in \operatorname{int} \left(\mathcal{E} \right) \quad \quad \Leftrightarrow \quad \quad \left. \eta \right|_{N} \ \text{is injective}. \label{2952}
\end{align}
\end{remark}

The functions in the evolutionary setting in this article typically depend on the time $t \in [0,\infty)$ in addition to a spatial variable $X \in \Omega_0$. We therefore make use of the concept of Bochner spaces, which generalize the classical $L^p$ spaces to spaces of functions the values of which lie in a general Banach space instead of in $\mathbb{R}$. More specifically, if $V(\Omega_0)$ is a Banach space of functions defined on $\Omega_0$ with a norm $\|\cdot \|_{V(\Omega_0)}$, then the Bochner space $L^p(0,T;V(\Omega_0))$, where $T \in (0,\infty]$ and $p \in [1,\infty]$, consists of (measurable) functions $m: [0,T] \rightarrow V(\Omega_0)$ for which the norm
\begin{align}
\left\| m \right\|_{L^p (0,T;V(\Omega_0))} := \left\{
\begin{array}{ll}
\left( \int_0^T \left\| m(t) \right\|_{V(\Omega_0)}^p \right)^\frac{1}{p} \quad \quad &\text{if } 1 \leq p < \infty, \\
\operatorname{ess\,sup}_{t \in [0,T]} \left\| m(t) \right\|_{V(\Omega_0)} \quad \quad &\text{if } p = \infty,\ 
\end{array}
\right. \label{bochnernorm}
\end{align}
is finite. This space itself, when equipped with the norm \eqref{bochnernorm}, constitutes a Banach space; for more details we refer e.g.\ to \cite[Section 1.5]{roubicek}. For our purposes, $V(\Omega_0)$ will typically be some type of Sobolev space on $\Omega_0$, meaning that we deal with functions that are $L^p$-integrable in time and weakly differentiable in space.

In order to be able to also work in the current configuration, we further require a generalization of these Bochner spaces to the setting of the moving domain $\Omega (\cdot)$. For the definition of such a generalization we assume the deformation $\eta$ to satisfy the conditions
\begin{align}
\eta \in L^\infty (0,T;\mathcal{E}), \quad \eta (t) \in \operatorname{int}\left( \mathcal{E} \right)\ \text{and} \ \tilde{E}_{\text{el}}\left( \eta (t) \right) \leq c \ \text{for a.a. } t \in [0,T], \label{2964}
\end{align}
where $c>0$ denotes a constant independent of $t \in [0,T]$ and $T \in (0,\infty]$ will be specified later. In particular, by the Ciarlet-Nečas condition \eqref{4016}, the mapping $X \mapsto \eta (t,X)$, $t \in [0,T]$, is injective in $\Omega_0$ and thus possesses an inverse
\begin{align}
\eta^{-1}(t,\cdot): \Omega \left(t \right) \rightarrow \Omega_0. \nonumber 
\end{align}

Then, for values $1 \leq p \leq \infty$, $1 \leq r < \infty$ and $k=0,1$, we define the generalized Bochner space
\begin{align}
L^p &\left(0 ,T;W^{k,r}\left(\Omega \left(\cdot \right)\right) \right) \nonumber \\
:=& \left\lbrace m : [0,T] \rightarrow \bigcup_{t \in [0,T]} W^{k,r}\left(\Omega(t)\right):\ m \left(\cdot, \eta^{-1}(\cdot, \cdot) \right) \in L^p \left(0,T;W^{k,r}\left(\Omega_0 \right) \right) \right\rbrace, \label{2961}
\end{align}

where the union is taken over uncountably many sets. Under the assumptions \eqref{2964} this space turns out to be a Banach space with the norm
\begin{align}
\left\| m \right\|_{L^p (0,T;W^{k,r}(\Omega (\cdot )))} := \left\{
\begin{array}{ll}
\left( \int_0^T \left\| m(t) \right\|_{W^{k,r}(\Omega (t))}^p \right)^\frac{1}{p} \quad \quad &\text{if } 1 \leq p < \infty, \\
\operatorname{esssup}_{t \in [0,T]} \left\| m(t) \right\|_{W^{k,r}(\Omega (t))} \quad \quad &\text{if } p = \infty,\ 
\end{array}
\right. \nonumber
\end{align}

as can be seen via a transformation to the reference configuration, cf.\ \cite[Lemma A.7.3]{thesis} and \cite[Theorem 2.4]{alphonse}.  Finally, for the weak formulation of the system \eqref{2949}--\eqref{2947}, presented in Section \ref{weaksolmagnetoel} below, we define the stray field associated to the magnetization via the variational formulation of the Poisson problem \eqref{2932}. To this end we introduce the notation
\begin{align}
\dot{H}^1\left( \mathbb{R}^3 \right) := \tilde{H}^1\left(\mathbb{R}^3 \right) / \mathbb{R}, \label{sobolevspacer3}
\end{align}
where 
\begin{align}
\tilde{H}^1\left(\mathbb{R}^3 \right) := \left\lbrace \psi \in H^1_{\text{loc}}\left( \mathbb{R}^3 \right):\ \nabla_x \psi \in L^2 \left(\mathbb{R}^3 \right) \right\rbrace \nonumber
\end{align}
represents the space of local $H^1$-functions the gradient of which is square integrable over the whole space $\mathbb{R}^3$, cf.\ \cite[Section 3]{praetorius}. The quotient space $\dot{H}^1( \mathbb{R}^3)$, in which the constant functions from the space $\tilde{H}^1(\mathbb{R}^3)$ are factored out, constitutes a Hilbert space with the bilinear product
\begin{align}
\left\langle \phi, \psi \right\rangle_{\dot{H}^1(\mathbb{R}^3) \times \dot{H}^1(\mathbb{R}^3)} := \int_{\mathbb{R}^3} \nabla_x \phi \cdot \nabla_x \psi \ dx, \nonumber
\end{align}
see \cite[Lemma 3.2]{praetorius}.The space $\dot{H}^1(\mathbb{R}^3)$ is also referred to as a homogeneous Sobolev space, for more details we refer to \cite[Section II.6]{galdi}. Now, for a given magnetization $\tilde{M} \in H^1(\Omega_0)$ in the reference configuration and a given deformation $\eta \in \mathcal{E}$ we denote the associated stray field by $H[\tilde{M}, \eta] = - \nabla_x \phi [\tilde{M}, \eta] \in L^2(\mathbb{R}^3)$, where $\phi[\tilde{M}, \eta] \in \dot{H}^1(\mathbb{R}^3) \bigcap H_{\text{loc}}^2(\eta( \Omega_0))$ satisfies the weak formulation
\begin{align}
\int_{\mathbb{R}^3} \nabla_x \phi\left[\tilde{M}, \eta \right] \cdot \nabla_x \psi \ dx = \int_{\eta (\Omega_0)} M \cdot \nabla_x \psi \ dx \quad \quad \forall \psi \in \dot{H}^1\left(\mathbb{R}^3 \right) \label{3008}
\end{align}

of the Poisson equation \eqref{2932}, cf.\ Lemma \ref{poissonequation}.

\section{Solution approach} \label{approximatesystemmagnetoelastic}

The main result of the present article, which is stated in Theorem \ref{mainresultmagnetoelastic} in the following section, is the proof of the existence of weak solutions to the model \eqref{2949}--\eqref{2947}. This proof, which constitutes a shortened version of the proof of the same result already presented in the PhD thesis \cite{thesis} of the third author, is carried out rigorously in Section \ref{proofsection} below. In the present section, we give an overview of the ideas of the proof and relate it to some physical interpretations.

Our proof is based on De Giorgi's minimizing movements scheme (see \cite{degiorgi}): We introduce a time-discrete approximation to the original problem by choosing a parameter $\Delta t > 0$ and establish a set of discrete times $k\Delta t$, $k \in \mathbb{N}_0$, partitioning the interval $[0,\infty)$ into intervals $[(k-1)\Delta t,k \Delta t)$. Then, at each discrete time $k\Delta t$ we solve a minimization problem and obtain discrete approximations of the equation of motion \eqref{2949} and the magnetic force balance \eqref{2948} as the corresponding Euler-Lagrange equations. Subsequently, passing to the limit in these equations, we recover a solution to the original system.

Before giving a detailed explanation of our approach we first present the full approximate problem: We fix an arbitrary discrete time $k \Delta t$, $k \in \mathbb{N}$. Beginning with the initial data, we assume that the solution $(\eta_{\Delta t}^{k-1}, \tilde{M}_{\Delta t}^{k-1}) \in \mathcal{E} \times H^1(\Omega_0)$ to the approximate problem at the time $(k-1)\Delta t$ has already been determined. We then consider the minimization problem
\begin{align}
\text{Find a minimizer} \quad \left( \eta_{\Delta t}^{k}, \tilde{M}_{\Delta t}^{k} \right) \in \mathcal{E} \times H^1\left(\Omega_0 \right) \quad \text{of} \quad \tilde{F}_{\Delta t}^k\left( \cdot, \cdot \right) \quad \text{over} \quad \mathcal{E} \times H^1\left(\Omega_0 \right), \label{3016}
\end{align}

where the functional $\tilde{F}_{\Delta t}^k: \mathcal{E} \times H^1(\Omega_0) \rightarrow \mathbb{R}$ is defined by
\begin{align}
\tilde{F}_{\Delta t}^k\left( \eta, \tilde{M} \right) :=& \tilde{E} \left( \eta, \tilde{M} \right) + \Delta t \tilde{R} \left( \eta_{\Delta t}^{k-1}, \frac{\eta- \eta_{\Delta t}^{k-1}}{\Delta t}, \frac{\tilde{M}- \tilde{M}_{\Delta t}^{k-1}}{\Delta t} \right) \nonumber \\
&- \int_{\Omega_0} \Delta t \rho f_{\Delta t}^k \left( \eta_{\Delta t}^{k-1} \right) \cdot \left( \frac{\eta - \eta_{\Delta t}^{k-1}}{\Delta t} \right) + \mu \tilde{M} \cdot \left(H_{\text{ext}}\right)_{\Delta t}^k(\eta)\ dX \label{2971}
\end{align}

with the energy potential $\tilde{E}$ defined in \eqref{2999} and the dissipation potential $\tilde{R}$ defined in \eqref{2936}. The stray field is discretized implicitly via the Poisson equation \eqref{3008} associated to $\tilde{M}$ and $\eta$. The discrete approximation $f_{\Delta t}^k$ of the given function $f$ at the time $k \Delta t$ in the formula \eqref{2971} is defined as
\begin{align}
f_{\Delta t}^k := f_{\kappa}(k\Delta t), \quad f_\kappa (t):= \int_0^T \theta _\kappa \left( t + \xi_\kappa (t)) - s \right)f(s)\ ds,\quad \xi_\kappa (t) := \kappa \frac{T-2t}{T}, \label{2966}
\end{align}

for a mollifier $\theta_\kappa :\mathbb{R}\to \mathbb{R}$ with support in $[-\kappa, \kappa]$ and a suitable choice of $\kappa = \kappa ( \Delta t),$ $\kappa (\Delta t)\rightarrow 0$ for $\Delta t \rightarrow 0$. The discrete approximation $(H_{\text{ext}})_{\Delta t}^k$ of the given function $H_{\text{ext}}$ at the time $k \Delta t$ instead is defined as the zero-order Clément quasi-interpolant
\begin{align}
(H_{\text{ext}})_{\Delta t}^k := \frac{1}{\Delta t} \int_{(k-1)\Delta t}^{k\Delta t} H_{\text{ext}} (s)\ ds, \label{2965}
\end{align}

cf.\ \cite[Remark 8.15]{roubicek}. \\

Next we discuss the mechanical and mathematical ideas behind the above approximate problem: The general procedure in our proof follows closely the implementation of De Giorgi's minimizing movements scheme (cf.\ \cite{degiorgi}) as used in \cite[Section 2]{variationalapproachtofsi}. This approach is what allows us to handle the coupling between the equation of motion \eqref{2949} and the magnetic force balance \eqref{2948}. Indeed, at each fixed discrete time an approximate deformation and an approximate magnetization are determined simultaneously by solving only one single minimization problem, namely the problem \eqref{3016}. From a physical point of view, this means that in our model magnetization and elastic deformation are coupled so that they reduce the energy together.

The existence of a minimizer to the problem \eqref{3016} can be shown via the direct method of the calculus of variations, cf.\ Lemma \ref{solutionminimization} below. Both a discrete approximation of the equation of motion and a discrete approximation of the magnetic force balance are then obtained as the Euler-Lagrange equations associated to the minimization problem; see the relations \eqref{3062} and \eqref{3034} below. In particular, this method allows us to evade the utilization of fixed point arguments, which are usually used to solve coupled systems of partial differential equations and which are typically designed for convex state spaces. In our setting, the latter condition is not satisfied since we consider deformations from the non-convex set $\mathcal{E}$ defined in \eqref{statespace}.

\begin{remark} \label{minimizingmovementsandnonconvexity}
De Giorgi's minimizing movements scheme moreover helps us to deal with the non-convexity of the energy: If we solved a discrete approximation of the equation of motion directly instead of via minimization, the classical way to obtain uniform (with respect to $\Delta t$) estimates for the deformation would be to test the equation at each time $k\Delta t$ by $\eta_{\Delta t}^k - \eta_{\Delta t}^{k-1}$. Then the desired bounds could be concluded provided that a discrete chain rule of the form
\begin{align}
\tilde{E} \left( \eta_{\Delta t}^k, \tilde{M} \right) - \tilde{E} \left( \eta_{\Delta t}^{k-1}, \tilde{M} \right) \leq \int_{\Omega_0} \tilde{E}_\eta \left( \eta_{\Delta t}^k, \tilde{M} \right) \cdot \left( \eta_{\Delta t}^k - \eta_{\Delta t}^{k-1} \right) \ dX \quad \quad \forall \tilde{M} \in H^1\left( \Omega_0 \right) \nonumber
\end{align}

holds true. However, due to the non-convexity of $\tilde{E}$ in the first argument (due to the term \eqref{2935}), such an estimate cannot be guaranteed. Since $\tilde{E}$ is further non-convex in the second argument (cf.\ the term \eqref{2934}), a corresponding problem also arises for the magnetization. Nonetheless, the De Giorgi method provides the necessary bounds for the discrete solution in a different way. More precisely, we obtain a uniform estimate for both the energy and the (discrete) dissipation potential by comparing the value of the functional $\tilde{F}_{\Delta t}^k$ in its minimizer $(\eta_{\Delta t}^k, \tilde{M}_{\Delta t}^k)$ to its value in the pair $(\eta_{\Delta t}^{k-1}, \tilde{M}_{\Delta t}^{k-1})$, cf.\ Lemma~\ref{uniformbounds} below.
\end{remark}

In the construction of the functional $\tilde{F}_{\Delta t}^k$ to the minimization problem \eqref{3016} it is crucial that its variation with respect to the deformation yields a suitable approximation of the equation of motion while its variation with respect to the magnetization leads to a suitable approximation of the magnetic force balance. Typically, this is achieved by composing the functional of the energy potential and a discretized version of the dissipation potential from the continuous problem. 

When selecting the discretized dissipation potential, it is important to note that, while the continuous equations are obtained by taking the variation of the dissipation with respect to time derivatives, the approximate equations on the time-discrete level can only be obtained by taking the variation with respect to the corresponding state variables themselves. A classical approach consists of choosing the discrete dissipation potential as the continuous dissipation potential evaluated in the state variables from the previous discrete time and the corresponding discrete time derivatives, represented by the differences between the current state variables and the state variables from the previous discrete time, divided by $\Delta t$. In our setting this means the choice
\begin{align}
\tilde{R} \left( \eta_{\Delta t}^{k-1}, \frac{\eta- \eta_{\Delta t}^{k-1}}{\Delta t}, \frac{\tilde{M}- \tilde{M}_{\Delta t}^{k-1}}{\Delta t} \right). \nonumber
\end{align}

The main difficulty in our problem lies in the fact that for the continuous system written in the form \eqref{2949}--\eqref{2947} it is not obvious that both equations can be expressed via the same energy and dissipation potentials. Indeed, for this to be possible, the transport \eqref{2938} of the magnetization in the magnetic force balance \eqref{2948} needs to appear in this dissipation potential without giving a contribution to the equation of motion. Hence, it is not clear after which energy and dissipation the functional $\tilde{F}_{\Delta t}^k$ should be modeled. This difficulty, however, can be tackled by the term
\begin{align}
\frac{1}{2} \left| \frac{1}{\det \left(\nabla \eta \right)} \partial_t \tilde{M} \right|^2 \det \left( \nabla \eta \right) \nonumber
\end{align}

in the dissipation potential $\tilde{R}$ in \eqref{2936}. This term, being independent of $\partial_t \eta$, vanishes when the variation of $\tilde{R}$ is taken with respect to $\partial_t \eta$. Therefore, it plays no role in the Piola-Kirchhoff stress tensor $\sigma$ (cf.\ \eqref{2933}) in the equation of motion \eqref{2949}. Yet, when $\tilde{R}$ is transformed to the current configuration (cf.\ \eqref{alternativedissipation}), the term turns into $\frac{1}{2}|D_tM|^2$. This allows us to express the magnetic force balance in the desired form \eqref{alternativemagforcebalance}, in which the energy potential $E$ and the dissipation potential $R$ coincide with the corresponding potentials $\tilde{E}$ and $\tilde{R}$ from the equation of motion \eqref{2949} formulated in the reference configuration. This knowledge offers us the opportunity to model the functional $\tilde{F}_{\Delta t}^k$ for our discrete minimization problem on the basis of $\tilde{E}$ and $\tilde{R}$, which is what we do in \eqref{2971}.

\begin{remark}
More intuitively, the discrete dissipation potential could also be expressed in the current configuration (at the time $(k-1)\Delta t$). The resulting more complicated expression
\begin{align}
& \int_{\Omega_{\Delta t}^{k-1}} \nu \left| \nabla_x \left( \frac{\eta \left( \left( \eta_{\Delta t}^{k-1} \right)^{-1}\right) - x}{\Delta t} \right) \right|^2 \nonumber \\
&+ \frac{1}{2} \left| \frac{1}{\det \left(\left[ \eta_{\Delta t}^{k-1} \right]\left( \eta_{\Delta t}^{k-1} \right)^{-1} \right)} \cdot \frac{\tilde{M}_{\Delta t}^k \left( \left( \eta_{\Delta t}^{k-1} \right)^{-1} \right) - \tilde{M}_{\Delta t}^{k-1} \left( \left( \eta_{\Delta t}^{k-1} \right)^{-1} \right)}{\Delta t} \right|^2 \ dx, \nonumber 
\end{align}
where $\Omega_{\Delta t}^{k-1} := \eta_{\Delta t}^{k-1}(\Omega_0)$, however, provides no added value to our proof. We point out that, due to the formulation of the minimization problem in Lagrangian coordinates, the Euler-Lagrange equation obtained by taking the variation of $\tilde{F}_{\Delta t}^k$ with respect to the magnetization is formulated in the reference configuration, despite the original magnetic force balance \eqref{2948} being formulated in the current configuration. This, however, is in fact favorable, as the limit passage with respect to $\Delta t \rightarrow 0$ is more convenient in Lagrangian coordinates. The final magnetic force balance is obtained afterwards by a simple transformation of the resulting limit identity to the current configuration.
\end{remark}

\begin{comment}
We remark that the Euler-Lagrange equation obtained by taking the variation of $\tilde{F}_{\Delta t}^k$ with respect to the magnetization is formulated in the reference configuration, despite the original magnetic force balance \eqref{3057} being formulated in the current configuration. This is in fact favorable, as the limit passage with respect to $\Delta t \rightarrow 0$ is more convenient in Lagrangian coordinates. The final magnetic force balance \eqref{3057} is obtained afterwards by a simple transformation of the resulting limit identity to the current configuration.
\end{comment}

Finally, we point out that the reason for discretizing $H_{\operatorname{ext}}$ via the zero-order Clément quasi interpolant \eqref{2965} lies in the derivation of the energy estimate for the discretized system in Lemma \ref{uniformbounds} below: In order to deduce a bound independent of $\Delta t$ we need to control the difference quotient of the chosen discretization of $H_{\text{ext}}$ through the classical time derivative of $H_{\text{ext}}$. For the choice \eqref{2965} of the discretization of $H_{\text{ext}}$ this is possible via Lemma \ref{differencequotientestimate} in the appendix.

\section{Weak solutions and main result} \label{weaksolmagnetoel}

In this section we introduce our definition of weak solutions to the system \eqref{2949}--\eqref{2947} and present our main result. First, however, we impose some additional conditions on the anisotropy energy density $\tilde{\Psi}$ and the elastic energy density $W$. More specifically, we assume that
\begin{align}
\tilde{\Psi}: \mathbb{R}^{3 \times 3} \times \mathbb{R}^3 \rightarrow \mathbb{R}_0^+,\quad W: \mathbb{R}^{3 \times 3} \rightarrow \mathbb{R}_0^+,&\quad \quad \tilde{\Psi} \in C^1 \left( \mathbb{R}^{3 \times 3} \times \mathbb{R}^3 \right),\quad W \in C^{1} \left( \mathbb{R}^{3 \times 3} \right), \label{3063} \\
W\left( \theta \right) &\geq c \left(\left| \theta \right|^{p_1} - 1 \right) \quad \quad \quad \quad \quad \forall \theta \in \mathbb{R}^{3\times 3}, \label{4003} \\
\left| W\left( \theta \right) \right| + \left| W'\left( \theta \right) \right| &\leq c \left( 1 + \left| \theta \right|^{p_2} \right)\quad \quad \quad \quad \quad \forall \theta \in \mathbb{R}^{3\times 3}, \label{3090} \\
\left| \tilde{\Psi} \left( \theta, \xi \right) \right| + \left| \tilde{\Psi}_F \left( \theta, \xi \right) \right| &\leq c \left( 1 + \left| \theta \right|^{p_2} + \left| \xi \right|^{p_3} \right) \quad \quad \forall \theta \in \mathbb{R}^{3\times 3},\ \xi \in \mathbb{R}^3, \label{2980} \\
\left| \tilde{\Psi}_M \left( \theta, \xi \right) \right| &\leq c \left( 1 + \left| \theta \right|^{p_2} + \left| \xi \right|^{p_4} \right) \quad \quad \forall \theta \in \mathbb{R}^{3\times 3},\ \xi \in \mathbb{R}^3 \label{3064}
\end{align}
for some constants $p_1,p_2,p_3,p_4 \in \mathbb{R}$ satisfying $2 \leq p_1 < \infty$, $1 \leq p_2 < \infty$, $1 \leq p_3 < 6$, $1 \leq p_4 < 5$. Here $\tilde{\Psi}_F$, $\tilde{\Psi}_M$ denote the derivatives of $\tilde{\Psi}$ with respect to the first and the second variable respectively and $c> 0$ denotes a constant independent of $\theta$ and $\xi$.

All these assumptions essentially set the growth of the energy functions at "infinity", i.e. they become relevant only at very large strains/values of the magnetization or its gradient. Thus, they mostly serve the mathematical purpose of setting a space of admissible states, since the growth controls which power of the magnetization, the gradient of the magnetization and/or deformation gradient should have a finite integral over the domain $\Omega_0$. Loosely speaking, this, in turn, reflects how much they can locally concentrate. Nevertheless, concentrations of the deformation gradient, for instance, are expected to trigger further inelastic processes, which are not captured in our model and that would not allow the deformation to grow to infinity (even locally). The prescribed growth may provide an approximation to this idea, as it again controls how fast the possible growth can be.\\

We now provide a weak formulation of equations  \eqref{2949}--\eqref{2947} which essentially means that we rewrite the functional derivatives via test functions and dualities (that can be understood as a generalization of the scalar product) as well as integrate in time. The dual space of a function space $X$ is denoted by $X^*$.

\begin{definition}
\label{weaksolutionsmagnetoelastic}
Let $\Omega_0 \subset \mathbb{R}^3$ be a bounded domain of class $C^{0,1}$. 
\begin{comment}
$\Omega_0, K \subset \mathbb{R}^3$ be bounded domains of class $C^{0,1}$.
\end{comment}
Let $N \subset \partial \Omega_0$, assume $P := \partial \Omega_0 \setminus N$ to have positive $2$-dimensional Hausdorff measure $\mathcal{H}^2(P)>0$ and let $\gamma : P \rightarrow \mathbb{R}^3$ be a given injective boundary deformation which can be extended into $\Omega_0$ as a deformation with finite elastic energy. Let $\rho, A, \beta, \nu, \mu>0$, $q>3$ and $a > \frac{3q}{q-3}$ and consider some data 
$$f \in L^\infty ((0,\infty)\times \mathbb{R}^3), \quad H_{\operatorname{ext}} \in W^{1,\frac{4}{3}}(0,\infty;W^{1,\frac{4}{3}}(\mathbb{R}^3)), \quad \eta_0 \in \operatorname{int} (\mathcal{E})\quad  \mbox{and} \quad \tilde{M}_0 \in H^1(\Omega_0)$$ such that $\tilde{E}(\eta_0, \tilde{M}_0)< \infty$, where $\tilde{E}$ is defined in \eqref{2999}, and $\eta_0$ is injective on $\partial \Omega_0$. Assume further that the energy densities $\tilde{\Psi}$ and $W$ satisfy the conditions \eqref{3063}--\eqref{3064}. Then the system \eqref{2949}--\eqref{2947} is said to admit a weak solution on $[0,T)$ for some $T>0$ if there exist functions
\begin{align}
\eta \in L^\infty \left(0,T;\mathcal{E} \right) ,\quad \quad \tilde{M} \in L^\infty \left( 0,T;H^1\left(\Omega_0 \right) \right) \label{2958}
\end{align}

with
\begin{align}
\partial_t \eta \in L^2 \left(0,T;H^1\left(\Omega_0 \right) \right) ,\quad \quad \partial_t \tilde{M} \in L^2 \left( \left(0,T\right) \times \Omega_0 \right), \label{2957}
\end{align}

such that the following statements hold true: 
\begin{itemize}
\item[(i)] The pair $(\eta, \tilde{M})$ satisfies
\begin{align}
& \int_0^T \left\langle \tilde{E}_\eta \left( \eta, \tilde{M} \right), \chi \right\rangle_{\left( W^{2,q} \left(\Omega_0 \right) \right)^* \times W^{2,q} \left(\Omega_0 \right)} + \left\langle \tilde{R}_{\partial _t \eta} \left(\eta, \partial_t \eta, \partial_t \tilde{M} \right), \chi \right\rangle_{\left( H^1 \left(\Omega_0 \right) \right)^* \times H^1 \left(\Omega_0 \right)} \ dt \nonumber \\
& - \int_0^T \int_{\Omega_0} \rho f \left(\eta \right) \cdot \chi + \mu \left[\left( \nabla \left( H_{\operatorname{ext}} \left(\eta \right)\right) \left(\nabla \eta \right)^{-1} \right)^T\tilde{M} \right] \cdot \chi \ dXdt = 0 \label{3075}
\end{align}

for all $\chi \in \mathcal{D}((0,T)\times \Omega_0)$, where $\mathcal{D}((0,T)\times \Omega_0)$ denotes the space of all smooth and compactly supported functions in $(0,T)\times \Omega_0$, as well as the initial conditions
\begin{align}
\eta(0) = \eta_0,\quad \quad \tilde{M}(0) = \tilde{M}_0 \label{3100}
\end{align}

in the sense that
\begin{align}
\lim_{t \rightarrow 0+} \left\| \eta(t) - \eta_0 \right\|_{C^1(\overline{\Omega_0})} = 0,\quad \quad \lim_{t \rightarrow 0+} \left\| \tilde{M}(t)- \tilde{M}_0 \right\|_{L^2(\Omega_0)} = 0. \label{3101}
\end{align}

\item[(ii)]  The pair
\begin{align}
\left(\eta, M \right),\quad \quad M:= M_\eta \left[ \tilde{M} \right] = \frac{1}{\det \left( \left[ \nabla \eta \right] \left( \eta^{-1} \right) \right)} \tilde{M} \left( \eta^{-1} \right) \in L^\infty \left(0,T;H^1 \left( \Omega (\cdot) \right) \right), \label{2951}
\end{align}

where $\Omega (t):= \eta (t,\Omega_0)$, satisfies
\begin{align}
& \int_0^T \left\langle E_M \left(\eta, M \right), \hat{M} \right\rangle_{\left( H^1 \left(\Omega (t) \right) \right)^* \times H^1 \left(\Omega (t) \right)} + \left\langle R_{D_t M} \left( \eta, v, D_t M \right), \hat{M} \right\rangle_{L^2(\Omega (t))\times L^2(\Omega (t))} \ dt \nonumber \\
&- \int_0^T \int_{\Omega (t)} \mu H_{\operatorname{ext}} \cdot \hat{M} \ dxdt = 0 \label{3057} 
\end{align}

for any test function $\hat{M} \in L^\infty (0,T;H^1(\Omega (\cdot)))$
with the velocity field $v$ defined in the current configuration via the relation
\begin{align}
v \left(t, \eta (t,X) \right) = \partial_t \eta (t,X)\quad \forall (t,X) \in [0,T] \times \Omega_0. \nonumber
\end{align}
\end{itemize}

\end{definition}

\begin{remark}
The initial condition for the magnetization in Definition \ref{weaksolutionsmagnetoelastic}, which is formulated in the reference configuration in \eqref{3100} and \eqref{3101}, can be expressed equivalently in the current configuration. More specifically, for
\begin{align}
M_0 := \frac{1}{\det \left( \left[ \nabla \eta_0 \right] \left( \eta_0^{-1} \right) \right)} \tilde{M}_0 \left( \eta_0^{-1} \right), \nonumber
\end{align}

a transformation between the reference configuration and the current configuration shows that the initial condition for $\tilde{M}$ in \eqref{3100} and \eqref{3101} is equivalent to the relation
\begin{align}
\lim_{t \rightarrow 0+} \left\| M(t)- \tilde{M}_0 \left( \eta_0 \left(\eta^{-1}(t) \right) \right) \right\|_{L^2(\Omega (t))} = 0. \nonumber
\end{align}

As it is more convenient, however, we choose the formulation in the reference configuration in Definition \ref{weaksolutionsmagnetoelastic}.
\end{remark}

%\subsection{Main result}

The main result of this chapter, which proves the existence of weak solutions to the system \eqref{2949}--\eqref{2947} as introduced in Definition \ref{weaksolutionsmagnetoelastic}, reads as follows:

\begin{theorem}
\label{mainresultmagnetoelastic}
Let the assumptions of Definition \ref{weaksolutionsmagnetoelastic} be satisfied. Assume further that 
the energy densities $\tilde{\Psi}$ and $W$ satisfy the conditions \eqref{3063}--\eqref{3064}. Then there exists a time $T' > 0$ such that the system \eqref{2949}--\eqref{2947} admits a weak solution $(\eta, \tilde{M})$ on $[0,T')$ in the sense of Definition~\ref{weaksolutionsmagnetoelastic}. Moreover, the time $T'$ can be chosen as the first time at which there occurs a self-contact of the material or the energy rises to infinity. %such that $T'=\infty$, or $\liminf_{t \rightarrow T'} \tilde{E}(\eta (t), \tilde{M}(t)) = \infty$, or $\eta (T') \in \partial \mathcal{E}$.
\end{theorem}

\begin{remark} \label{compactsupptestfunctions}
Theorem \ref{mainresultmagnetoelastic} needs to be understood as a local result in the sense that the test functions $\chi$ in the weak formulation \eqref{3075} of the equation of motion are compactly supported. The reason for this lies in the regularity of the stray field: Indeed, the gradient of $H[\tilde{M},\eta]$ seems to be at best locally integrable in space, cf.\ Lemma \ref{poissonequation} in the appendix. For global integrability at least $C^2$-regularity of the current configuration $\Omega (t)$, $t \in [0,T]$, appears to be required. However, even if we assumed $\Omega_0$ to be of class $C^2$, such a regularity could not be expected since the deformation $\eta(t)$ as a $W^{2,q}(\Omega_0)$-function (cf.\ the elastic energy \eqref{4039}) is only known to possess $C^{1,\alpha}$-regularity, $\alpha = 1 - \frac{3}{q}$, in the spatial variable. A generalization of Theorem \ref{mainresultmagnetoelastic} to non-compactly supported test functions in the equation of motion remains an open problem for future research.
\end{remark}

\begin{comment}
\begin{remark}
The ending criterion $\eta (T') \in \partial \mathcal{E}$ in Theorem \ref{mainresultmagnetoelastic} can be understood as the occurrence of a self-contact of the material at the time $T'$, i.e.\ we prove the existence of a weak solution only up to the time of the first self-collision. The problem that occurs when $\eta (T') \in \partial \mathcal{E}$ is that the variation of the functional of the minimization problem, which we use to prove Theorem \ref{mainresultmagnetoelastic} (cf.\ Section \ref{approximatesystemmagnetoelastic}), cannot be taken in all directions anymore and thus it is not possible to calculate the full Euler-Lagrange equations. Nonetheless, there exist methods to prove the existence of weak solutions even beyond self-collision, as has been done for example for the evolution of viscoelastic solids in \cite{antonin1,antonin2,kromerroubicek}.
\end{remark}
\end{comment}

\begin{remark} \label{selfcontact}
Mathematically, the ending criterion of $T'$ being the time of the first self-contact in Theorem \ref{mainresultmagnetoelastic} means that $\eta (T') \in \partial \mathcal{E}$. The problem that occurs in this case is that the variation of the functional of the minimization problem, which we use to prove Theorem \ref{mainresultmagnetoelastic} (cf.\ Section \ref{approximatesystemmagnetoelastic}), cannot be taken in all directions anymore and thus it is not possible to calculate the full Euler-Lagrange equations. Nonetheless, there exist methods to prove the existence of weak solutions even beyond self-collision, as has been done for example for the evolution of viscoelastic solids in \cite{antonin1,antonin2,kromerroubicek}. The alternative ending criterion of the energy blowing up at the time $T'$ can be expressed as $\liminf_{t \rightarrow T_{\operatorname{max}}}\tilde{E}(\eta (t),\tilde{M}(t))< \infty.$ We point out that this possibility could be avoided through the addition of suitable stabilizing terms to the discrete approximation scheme as in \cite{mielkeroubicek}. If there occurs no self-contact and the energy remains finite at any finite time, the ending time $T'$ can further be chosen as $T' = \infty$.
\end{remark}

The proof of Theorem \ref{mainresultmagnetoelastic} will be achieved via an approximation method, which we present in the following section and which is carried out in Section \ref{proofsection}.

\section{Proof of the main result} \label{proofsection}

The remainder of the article is dedicated to the proof of our main result Theorem \ref{mainresultmagnetoelastic}. A more detailed version of this proof can be found in the PhD thesis \cite{thesis} of the third author. We begin by proving the existence of a solution to the approximate problem introduced in Section \ref{approximatesystemmagnetoelastic}.

\begin{proposition}
\label{solutionminimization}
Let all the assumptions of Theorem \ref{mainresultmagnetoelastic} be satisfied and let $\Delta t > 0$. Let further $f_{\Delta t}^k$ be given by \eqref{2966} and $(H_{\text{ext}})_{\Delta t}^k$ be given by \eqref{2965} for any $k \in \mathbb{N}_0$. Then, for all $k \in \mathbb{N}$, there exists a solution
\begin{align}
\left( \eta_{\Delta t}^k, \tilde{M}_{\Delta t}^k \right) \in \mathcal{E} \times H^1\left(\Omega_0\right) \nonumber
\end{align}

to the minimization problem \eqref{3016}.% Moreover, if $T>0$ is such that $\frac{T}{\Delta t} \in \mathbb{N}$ and $\eta_{\Delta t}^k \notin \partial \mathcal{E}$ for all $k=1,...,\frac{T}{\Delta t}$, then the approximate equation of motion \eqref{3062} is satisfied for all $\chi \in \mathcal{D}((0,T)\times \Omega_0)$ and the approximate magnetic force balance \eqref{3034} is satisfied for all $\tilde{\hat{M}} \in L^\infty (0,T;H^1(\Omega_0))$.
\end{proposition}

\textbf{Proof}

We fix some arbitrary discrete time index $k \in \mathbb{N}$ and assume that the minimization problem \eqref{3016} has been solved for each time index $l=1,...,k-1$. We argue via the direct method and thus first show that $\tilde{F}_{\Delta t}^k$ is bounded from below on $\mathcal{E}\times H^1(\Omega_0)$: For the stray field term appearing in the energy $\tilde{E}$ and thus in $\tilde{F}_{\Delta t}^k$, we note that, by the definition of $H[\tilde{M}, \eta]$ via the Poisson equation \eqref{3008},
\begin{align}
\int_{\Omega_0} - \frac{\mu}{2} \tilde{M} \cdot H \left[ \tilde{M}, \eta \right](\eta) \ dX = \int_{\eta(\Omega_0)} - \frac{\mu}{2} M_{\eta}\left[ \tilde{M} \right] \cdot H\left[ \tilde{M}, \eta \right]\ dx = \int_{\mathbb{R}^3} \frac{\mu}{2} \left| H \left[ \tilde{M}, \eta \right] \right|^2\ dx \geq 0. \label{3003}
\end{align}

Continuing with the $f$-dependent term in $\tilde{F}_{\Delta t}^k$, we remark that, by Lemma \ref{boundawayfromzero} in the appendix, the quantity $\det (\nabla \eta_{\Delta t}^{k-1})$ is bounded away from zero in dependence of only the value $\tilde{E}_{\operatorname{el}}(\eta_{\Delta t}^{k-1})$. This, in combination with Young's inequality and the Poincaré inequality, allows us to estimate
\begin{align}
&\int_{\Omega_0} \Delta t \nu \left| \nabla \left( \frac{\eta - \eta_{\Delta t}^{k-1}}{\Delta t} \right) \left( \nabla \eta_{\Delta t}^{k-1} \right)^{-1} \right|^2\det \left( \nabla \eta_{\Delta t}^{k-1} \right) - \Delta t \rho f_{\Delta t}^k \left( \eta_{\Delta t}^{k-1} \right) \cdot \left( \frac{\eta - \eta_{\Delta t}^{k-1}}{\Delta t} \right)\ dX \nonumber \\
\geq& \Delta t \left[\frac{\nu c\left( \eta_{\Delta t}^{k-1} \right)}{2} \int_{\Omega_0} \left| \nabla \left( \frac{\eta - \eta_{\Delta t}^{k-1}}{\Delta t}\right) \right|^2\ dX - \frac{\left| \Omega_0 \right| \left( \rho \left\| f \right\|_{L^\infty ((0,\infty)\times \mathbb{R}^3)} \right)^2 }{2\nu c\left( \eta_{\Delta t}^{k-1} \right)}\right] \geq -c(\Delta t, \eta_{\Delta t}^{k-1},f) \label{3005}
\end{align}

for a constant $c(\eta_{\Delta t}^{k-1})>0$ which remains bounded away from $0$ for bounded values of $\tilde{E}_{\text{el}}(\eta_{\Delta t}^{k-1})$ and a constant $c(\Delta t, \eta_{\Delta t}^{k-1},f) > 0 $, both independent of $\eta$ and $\tilde{M}$. In order to control the $H_{\text{ext}}$-dependent term in $\tilde{F}_{\Delta t}^k$, we first deduce from the Gagliardo-Nirenberg inequality, the Poincaré inequality and the Morrey embedding $W^{2,q}(\Omega_0) \subset W^{1,\infty} (\Omega_0)$ (since, by Definition \ref{weaksolutionsmagnetoelastic}, $q>3$, cf.\ \cite[Section 1.3.5.8]{novotnystraskraba}) that
\begin{align}
c_1 \left\| \eta \right\|^2_{W^{1,\infty}(\Omega_0)} \leq \int_{\Omega_0} \frac{1}{2} W \left( \nabla \eta \right) + \frac{1}{2q} \left| \nabla^2 \eta \right|^q\ dX + 1 \label{3000}
\end{align} %see iphone photo 15.08.2023

for a constant $c_1>0$. Furthermore, we observe that 
\begin{align}
\int_{\eta (\Omega_0)} \frac{1}{8\beta^2}\ dx = \int_{\Omega_0} \frac{\operatorname{det} \left( \nabla \eta \right)}{8 \beta^2}\ dX \leq c_1 \left\| \eta \right\|_{W^{1,\infty}(\Omega_0)}^2 + c \nonumber
\end{align}

for some constant $c>0$. This, together with \eqref{3000} and Young's inequality, allows us to absorb the $H_{\text{ext}}$-dependent term into the non-negative energy,
\begin{align}
&\int_{\Omega_0} \frac{1}{2} W \left( \nabla \eta \right) + \frac{1}{2q} \left| \nabla^2 \eta \right|^q + \frac{1}{8\beta^2} \left( \left| \frac{1}{\det \left( \nabla \eta \right)} \tilde{M} \right|^2 - 1 \right)^2 \det \left( \nabla \eta \right) - \mu \tilde{M} \cdot \left(H_{\text{ext}}\right)_{\Delta t}^k \left( \eta \right)\ dX \nonumber \\
\geq& c_1 \left\| \eta \right\|^2_{W^{1,\infty}(\Omega_0)} - 1 + \int_{\eta (\Omega_0)} \frac{1}{16\beta^2} \left| M_\eta \left[ \tilde{M} \right] \right|^4 - \frac{1}{8\beta^2} - \frac{1}{16\beta^2} \left| M_\eta \left[ \tilde{M} \right] \right|^4 - \frac{3\left(4\beta^2 \mu^4\right)^\frac{1}{3} \left| \left(H_{\text{ext}}\right)_{\Delta t}^k \right|^\frac{4}{3}}{4}\ dx \nonumber \\
\geq&  c_1 \left\| \eta \right\|^2_{W^{1,\infty}(\Omega_0)} - 1 -  c_1 \left\| \eta \right\|^2_{W^{1,\infty}(\Omega_0)} - c\left(H_{\text{ext}} \right) \geq - c\left(H_{\text{ext}} \right) \label{6}
\end{align}

for a constant $c( H_{\text{ext}})>0$ independent of $\eta$ and $\tilde{M}$.  The estimate \eqref{3003} shows that $\tilde{E}(\eta, \tilde{M})$ is non-negative. Hence, using the estimates \eqref{3005} and \eqref{6} to control the forcing terms in $\tilde{F}_{\Delta t}^k$ (cf.\ \eqref{2971}) as well as the definition of the energy $\tilde{E}$ in \eqref{2999}, we conclude that
\begin{align}
\tilde{F}_{\Delta t}^k \left( \eta, \tilde{M} \right) \geq \frac{1}{2} \tilde{E} \left( \eta, \tilde{M} \right) - c \geq -c > - \infty \label{4040}
\end{align}

for all $(\eta, \tilde{M}) \in \mathcal{E}\times H^1(\Omega_0)$ and a constant $c>0$ independent of $\eta$ and $\tilde{M}$. Consequently, there exists a minimizing sequence $(\eta_j, \tilde{M}_j)_{j \in \mathbb{N}}$ for $\tilde{F}_{\Delta t}^k$. Thanks to the estimate \eqref{3076} for the solution to the Poisson equation \eqref{3008}, it further holds that
\begin{align}
\left\| H \left[ \tilde{M}_j,\eta_j \right] \right\|_{L^2(\mathbb{R}^3)} = \left\| \phi \left[ \tilde{M}_j,\eta_j \right] \right\|_{\dot{H}^1(\mathbb{R}^3)} \leq c \nonumber
\end{align}
for a constant $c>0$ independent of $j$. This, together with the upper bound \eqref{4040} for the energy $\tilde{E}$, allows us to choose the minimizing sequence $(\eta_j, \tilde{M}_j)_{j \in \mathbb{N}}$ for $\tilde{F}_{\Delta t}^k$ in such a way that
\begin{align}
\eta_j \rightharpoonup \eta \quad \text{in } W^{2,q}\left(\Omega_0\right),\quad \quad \tilde{M}_j \rightharpoonup \tilde{M} \quad \text{in } H^1\left(\Omega_0\right),\quad \quad H\left[\tilde{M}_j, \eta_j \right] \rightharpoonup H \left[ \tilde{M}, \eta \right] \quad \text{in } L^2\left(\mathbb{R}^3 \right) \label{3012}
\end{align}

for some functions $\eta \in \mathcal{E}$, $\tilde{M} \in H^1(\Omega_0)$ and $H[\tilde{M},\eta] = - \nabla_x \phi[\tilde{M}, \eta]$, where $\phi[\tilde{M}, \eta]$ denotes the solution to the Poisson problem \eqref{3008} associated to $M$ and $\eta$. The latter identification could be achieved by passing to the limit in the corresponding Poisson equation satisfied by $\phi [\tilde{M}_j, \eta_j]$ under exploitation of the uniform convergence of $\eta_j$ implied by its $W^{2,q}$-convergence in \eqref{3012}. From the convergences \eqref{3012} and the continuity and boundedness assumptions \eqref{3063}, \eqref{3090} and \eqref{2980} on $W$ and $\tilde{\Psi}$ we infer that
\begin{align}
\tilde{F}_{\Delta t}^k\left( \eta, \tilde{M} \right) \leq \liminf_{j \rightarrow \infty} \tilde{F}_{\Delta t}^k \left( \eta_j, \tilde{M}_j \right). \label{3009}
\end{align}

This proves the existence of the desired solution $(\eta_{\Delta t}^k, \tilde{M}_{\Delta t}^k) = (\eta, \tilde{M})$ to the minimization problem \eqref{3016} at the discrete time $k\Delta t$ and thus the desired result.
$\hfill \Box$\\

Having proved the existence of a solution to the minimization problem \eqref{3016}, we are now interested in the Euler-Lagrange equations of the associated functional $\tilde{F}_{\Delta t}^k$ (cf.\ \eqref{2971}), which serve as discrete approximations of the equation of motion \eqref{3075} and the magnetic force balance \eqref{3057}. We assemble these equations as a time dependent system on some time interval $[0,T]$, where $T>0$ will be specified in Lemma \ref{injectivity} below. To this end we define piecewise affine and piecewise constant interpolants of the discrete quantities: For all time-independent functions $h_{\Delta t}^k$, $k \in \mathbb{N}_0$, we set
\begin{align}
h_{\Delta t}(t) &:= \left( \frac{t}{\Delta t} - (k-1)\right)h_{\Delta t}^k + \left(k - \frac{t}{\Delta t}\right) h_{\Delta t}^{k-1}\ \ \ &\forall& t \in ((k-1)\Delta t, k\Delta t],\quad k \in \mathbb{N}, \label{2931} \\
\overline{h}_{\Delta t}(t) &:= h_{\Delta t}^k \ \ \ &\forall& t \in ((k-1)\Delta t, k\Delta t],\quad k \in \mathbb{N}_0, \label{2930} \\
\overline{h}'_{\Delta t}(t) &:= h_{\Delta t}^{k-1} \ \ \ &\forall& t \in ((k-1)\Delta t, k\Delta t],\quad k\in \mathbb{N}. \label{2929}
\end{align}
For the discrete approximation of the equation of motion, we consider $T>0$ such that $\frac{T}{\Delta t} \in \mathbb{N}$ and $\eta_{\Delta t}^k \notin \partial \mathcal{E}$ for all $k=1,...,\frac{T}{\Delta t}$. The existence of such a time $T$ independent of $\Delta t$ is shown below in Lemma \ref{injectivity}. As a consequence, for all $k=1,...,\frac{T}{\Delta t}$, it holds that $\eta_{\Delta t}^k + \epsilon \chi(t) \in \mathcal{E}$ for any $t \in [0,T]$, $\chi \in \mathcal{D}((0,T) \times \Omega_0)$ and any sufficiently small $\epsilon > 0$. This allows us to take the variation of $\tilde{F}_{\Delta t}^k$ at its minimizer $(\eta_{\Delta t}^k, \tilde{M}_{\Delta t}^k)$ with respect to the deformation. We integrate the resulting Euler-Lagrange equation over $[(k-1)\Delta t, k\Delta t]$, sum over $k$ and obtain the approximate equation of motion
\begin{align}
& \int_0^T \left\langle D_1 \tilde{E} \left( \overline{\eta}_{\Delta t}, \overline{\tilde{M}}_{\Delta t} \right), \chi \right\rangle_{\left( W^{2,q} \left(\Omega_0 \right) \right)^* \times W^{2,q} \left(\Omega_0 \right)} \nonumber \\
+& \left\langle \Delta t D_1 \tilde{R}\left(\overline{\eta}'_{\Delta t}, \partial_t \eta_{\Delta t}, \partial_t\tilde{M}_{\Delta t} \right), \chi \right\rangle_{\left( H^1 \left(\Omega_0 \right) \right)^* \times H^1 \left(\Omega_0 \right)} \ dt - \int_0^T \int_{\Omega_0} \rho \overline{f}_{\Delta t} \left(\overline{\eta}'_{\Delta t} \right) \cdot \chi \ dxdt  = 0 \label{3062}
\end{align}
for all $\chi \in \mathcal{D}((0,T)\times \Omega_0)$, where $D_1 \tilde{E}(\eta, \tilde{M})$ and $D_1 \tilde{R}(\eta, \tilde{M})$ denote the Fréchet derivatives of $\tilde{E}$ and $\tilde{R}$, respectively, with respect to the first argument at $(\eta, \tilde{M})$.

For the derivation of the approximate magnetic force balance we take the variation of $\tilde{F}_{\Delta t}^k$, $k=1,...,\frac{T}{\Delta t}$, at $(\eta_{\Delta t}^k, \tilde{M}_{\Delta t}^k)$ with respect to the magnetization. Integrating the resulting Euler-Lagrange equation over $[(k-1)\Delta t, k\Delta t]$ and summing over $k$, we deduce the approximate magnetic force balance
\begin{align}
& \int_0^T \left\langle D_2 \tilde{E} \left(\overline{\eta}_{\Delta t}, \overline{\tilde{M}}_{\Delta t} \right), \tilde{\hat{M}} \right\rangle_{\left( H^1 \left(\Omega_0 \right) \right)^* \times H^1 \left(\Omega_0 \right)} + \left\langle \Delta t D_2 \tilde{R}\left(\overline{\eta}'_{\Delta t}, \partial_t \eta_{\Delta t}, \partial_t\tilde{M}_{\Delta t} \right), \tilde{\hat{M}} \right\rangle_{L^2(\Omega_0)\times L^2(\Omega_0)} \ dt = 0 \label{3034} 
\end{align}
for all $\tilde{\hat{M}} \in L^\infty (0,T;H^1(\Omega_0))$, where $D_2 \tilde{E}(\eta, \tilde{M})$ and $D_2 \tilde{R}(\eta, \tilde{M})$ denote the Fréchet derivatives of $\tilde{E}$ and $\tilde{R}$, respectively, with respect to the second argument at $(\eta, \tilde{M})$. \\

Next, in order to return from the discrete to the continuous setting and recover a weak solution to the original problem, we need to pass to the limit with respect to $\Delta t \rightarrow 0$ in the equations \eqref{3062} and \eqref{3034}. To this end, we first establish an interval $[0,T]$ with some $T>0$ independent of $\Delta t$ on which we are able to find an energy estimate for the discrete solution, uniform with respect to $\Delta t > 0$. More precisely, we show the following lemma.
\begin{lemma}
\label{uniformbounds}
There exists a time $T>0$, independent of $\Delta t$, and a constant $c>0$, independent of $\Delta t$ and $k=1,...,\frac{T}{\Delta t}$, such that
\begin{align}
\tilde{E} \left( \eta_{\Delta t}^k, \tilde{M}_{\Delta t}^k \right) + \Delta t \sum_{l=1}^k \tilde{R} \left( \eta_{\Delta t}^{l-1}, \frac{\eta_{\Delta t}^l - \eta_{\Delta t}^{l-1}}{\Delta t}, \frac{\tilde{M}_{\Delta t}^l- \tilde{M}_{\Delta t}^{l-1}}{\Delta t} \right) \leq c \quad \quad \forall k = 1,...,\frac{T}{\Delta t}. \label{3029}
\end{align}
\end{lemma}

\begin{comment}
The proof of this lemma is mostly standard with the only difficulties being caused by the term depending on the external forcing term $f$ in the discrete functional $\tilde{F}_{\Delta t}^k$. Indeed, in order to control this term during the derivation of the estimate \eqref{3029} for some fixed $k$, we already need to know a uniform bound of $\det (\nabla \eta_{\Delta t}^{k-1})$ away from zero, cf.\ the deduction of the estimate \eqref{3005}. We achieve this via an induction argument, allowing us to assume $\tilde{E}_{\operatorname{el}}(\eta_{\Delta t}^{k-1})$ to be bounded.\\
\end{comment}

\textbf{Proof}

We choose some number $E_0 > 0$ such that
\begin{align}
\tilde{E} (\eta_0, \tilde{M}_0) < c_2 < E_0, \nonumber
\end{align}

for the constant $c_2>0$ chosen below in \eqref{97}, dependent on the data but independent of $\Delta t >0$. We further choose a time $T=T(E_0)>0$ sufficiently small such that
\begin{align}
\left[c_2 + T \frac{ \left| \Omega_0 \right| \left( \rho \left\| f \right\|_{L^\infty ((0,\infty)\times\mathbb{R}^3)} \right)^2 }{2\nu c_3} \right]e^{\max\{1,\ 8\beta^2\} T} \leq E_0, \label{3022}
\end{align}

where $c_3 > 0$ denotes the constant chosen below in \eqref{2990}, dependent on $E_0$ and the data but independent of $\Delta t >0$ and $k=1,...,\frac{T}{\Delta t}$. We argue via induction with respect to the discrete time index $k$. More specifically, we choose an arbitrary discrete time index $k=1,...,\frac{T}{\Delta t}$ and assume that
\begin{align}
\frac{1}{2}\tilde{E} \left( \eta_{\Delta t}^{l}, \tilde{M}_{\Delta t}^{l} \right)\leq E_0 \quad \quad \forall l = 0,...,k-1. \label{3021}
\end{align}

Each pair $(\eta_{\Delta t}^l, \tilde{M}_{\Delta t}^l)$, $l=1,...,k$, as a minimizer of the functional $\tilde{F}_{\Delta t}^l$, satisfies $\tilde{F}_{\Delta t}^l(\eta_{\Delta t}^l, \tilde{M}_{\Delta t}^l) \leq \tilde{F}_{\Delta t}^l(\eta_{\Delta t}^{l-1}, \tilde{M}_{\Delta t}^{l-1})$. We sum this inequality over all indices $l=1,...,k$ to infer from the definition of $\tilde{F}_{\Delta t}^l$ in \eqref{2971} that
\begin{align}
&\tilde{E}\left( \eta_{\Delta t}^k, \tilde{M}_{\Delta t}^k \right) + \Delta t \sum_{l=1}^k\tilde{R} \left( \eta_{\Delta t}^{l-1}, \frac{\eta_{\Delta t}^l - \eta_{\Delta t}^{l-1}}{\Delta t}, \frac{\tilde{M}_{\Delta t}^l- \tilde{M}_{\Delta t}^{l-1}}{\Delta t} \right) \nonumber \\
\leq& \tilde{E} \left(\eta_0, \tilde{M}_0\right) + \sum_{l=1}^k \left[ \int_{\Omega_0} \mu \tilde{M}_{\Delta t}^l \cdot \left(H_{\text{ext}}\right)_{\Delta t}^l\left(\eta_{\Delta t}^l\right)\ dX - \int_{\Omega_0} \mu \tilde{M}_{\Delta t}^{l-1} \cdot \left(H_{\text{ext}}\right)_{\Delta t}^l\left(\eta_{\Delta t}^{l-1}\right)\ dX \right] \nonumber \\
&+ \Delta t \sum_{l=1}^k \int_{\Omega_0} \rho f_{\Delta t}^l \left( \eta_{\Delta t}^{l-1} \right) \cdot \left( \frac{\eta_{\Delta t}^l - \eta_{\Delta t}^{l-1}}{\Delta t} \right) \ dX. \label{3095}
\end{align}

The $H_{\operatorname{ext}}$-dependent terms on the right-hand side of this inequality can be controlled via summation by parts, the bound of the discrete difference quotient of $H_{\operatorname{ext}}$ in terms of $\partial_t H_{\operatorname{ext}}$ (cf.\ Lemma \ref{differencequotientestimate} in the appendix) and the bound \eqref{3000} of the $W^{1,\infty}(\Omega_0)$-norm of the deformation in terms of the elastic energy:
\begin{align}
&\frac{1}{2}\tilde{E} \left( \eta_{\Delta t}^k, \tilde{M}_{\Delta t}^k \right) + \Delta t \sum_{l=1}^k \tilde{R} \left( \eta_{\Delta t}^{l-1}, \frac{\eta_{\Delta t}^l - \eta_{\Delta t}^{l-1}}{\Delta t}, \frac{\tilde{M}_{\Delta t}^l- \tilde{M}_{\Delta t}^{l-1}}{\Delta t} \right) \nonumber \\
\leq& c_2 + \Delta t \sum_{l=1}^{k-1} \max \left\lbrace 1,\ 8 \beta^2 \right\rbrace \frac{1}{2} \tilde{E} \left( \eta_{\Delta t}^l, \tilde{M}_{\Delta t}^l \right) + \Delta t \sum_{l=1}^k \int_{\Omega_0} \rho f_{\Delta t}^l \left( \eta_{\Delta t}^{l-1} \right) \cdot \left( \frac{\eta_{\Delta t}^l - \eta_{\Delta t}^{l-1}}{\Delta t} \right) \ dX, \label{2946}
\end{align}

where the constant
\begin{align}
c_2 = c_2\left(\eta_0, \tilde{M}_0, \beta, \mu, \Omega_0, q, p_1, \gamma, \left\| H_{\text{ext}} \right\|_{W^{1,\frac{4}{3}}(0,\infty;L^\frac{4}{3}(\mathbb{R}^3))} \right) > 0 \label{97}
\end{align}

is independent of $\Delta t$ and $k$. Moreover, the $f$-dependent terms on the right-hand side of this inequality can be controlled by recalling the first inequality in \eqref{3005}: Indeed, the constant $c(\eta_{\Delta t}^{k-1})$ in this estimate is bounded away from zero for bounded values of $\tilde{E}\left( \eta_{\Delta t}^{k-1}, \tilde{M}_{\Delta t}^{k-1} \right)$. In our current situation, due to the induction assumption \eqref{3021}, we can replace this constant by a constant
\begin{align}
c_3 = c_3(E_0, \eta_0, \tilde{M}_0, \mu, \beta, H_{\text{ext}})>0 \label{2990}
\end{align}

independent of $\Delta t$ and $k$. It follows that
\begin{align}
&\frac{1}{2}\tilde{E} \left( \eta_{\Delta t}^k, \tilde{M}_{\Delta t}^k \right) \nonumber \\
& + \Delta t \sum_{l=1}^k \left[ \frac{\nu c_3}{2} \int_{\Omega_0} \left| \nabla \left(\frac{\eta_{\Delta t}^l - \eta_{\Delta t}^{l-1}}{\Delta t}\right) \right|^2 + \det \left( \nabla \eta _{\Delta t}^{l-1} \right) \frac{1}{2} \left| \frac{1}{\det \left( \nabla \eta_{\Delta t}^{l-1}\right)} \cdot \frac{\tilde{M}_{\Delta t}^l - \tilde{M}_{\Delta t}^{l-1}}{\Delta t} \right|^2\ dX \right] \nonumber \\
\leq& c_2 + k\Delta t \frac{ \left| \Omega_0 \right| \left( \rho \left\| f \right\|_{L^\infty ((0,\infty)\times\mathbb{R}^3)} \right)^2 }{2\nu c_3} + \Delta t \sum_{l=1}^{k-1} \max \left\lbrace 1,\ 8 \beta^2 \right\rbrace \frac{1}{2} \tilde{E} \left( \eta_{\Delta t}^l, \tilde{M}_{\Delta t}^l \right). \nonumber
\end{align}

Now, the discrete Gronwall estimate (cf.\ \cite[(1.67)]{roubicek}) shows that the left-hand side of this inequality is bounded by
\begin{align}
\left[c_2 + k \Delta t \frac{ \left| \Omega_0 \right| \left( \rho \left\| f \right\|_{L^\infty ((0,\infty)\times\mathbb{R}^3)} \right)^2 }{2\nu c_3} \right]e^{\max \{1,\ 8\beta^2\} k\Delta t} \leq E_0, \nonumber
\end{align}

where the last inequality is due to our choice \eqref{3022} of $T$ and the fact that $k\Delta t \leq T$. This concludes the proof.
$\hfill \Box$ \\

For the time $T>0$ given by Lemma \ref{uniformbounds} the energy estimate \eqref{3029} yields several uniform bounds, which we summarize in the following corollary:

\begin{corollary} \label{boundscor}
Let $T>0$ be given by Lemma \ref{uniformbounds}. There exists a constant $c>0$, independent of $\Delta t>0$, such that for all $k=1,...,\frac{T}{\Delta t}$ it holds that
\begin{align}
\left\| \eta_{\Delta t}^k \right\|_{W^{2,q}(\Omega_0)} + \left\| \frac{1}{\det \left( \nabla \eta_{\Delta t}^k \right)} \right\|_{L^\infty (\Omega_0)} + \left\| \tilde{M}_{\Delta t}^k \right\|_{H^1(\Omega_0)} + \left\| \frac{\tilde{M}_{\Delta t}^k}{\det \left( \nabla \eta_{\Delta t}^k \right)} \right\|_{H^1(\Omega_0)} \leq c \label{uniformestimates}
\end{align}

and
\begin{align}
\left\| \partial_t \eta_{\Delta t} \right\|_{L^2(0,T;H^1(\Omega_0))} + \left\| \partial_t \tilde{M}_{\Delta t} \right\|_{L^2((0,T) \times \Omega_0)} &\leq c. \label{3035}
\end{align}
\end{corollary}

For the interpolated functions defined via the formulas \eqref{2931}--\eqref{2929}, the bounds given by Corollary \ref{boundscor} and the Aubin-Lions Lemma imply the existence of functions $\eta \in L^\infty(0,T;\mathcal{E})$, $\tilde{M} \in L^\infty(0,T;H^1(\Omega_0))$ such that, possibly after the extraction of a subsequence,
\begin{align}
\overline{\eta}_{\Delta t}, \overline{\eta}'_{\Delta t} \buildrel\ast\over\rightharpoonup \eta \ \text{in } L^\infty \left(0,T;W^{2,q}\left(\Omega_0\right)\right), \ \quad \quad \quad \quad \quad \partial_t \eta_{\Delta t} \rightharpoonup \partial_t \eta \ &\text{in } L^2\left(0,T;H^1\left(\Omega_0\right) \right), \label{40} \\
\overline{\tilde{M}}_{\Delta t}, \overline{\tilde{M}}'_{\Delta t} \buildrel\ast\over\rightharpoonup \tilde{M} \ \text{in } L^\infty \left(0,T;H^1 \left(\Omega_0\right)\right),\quad \quad \quad \quad \ \ \ 
\partial_t \tilde{M}_{\Delta t} \rightharpoonup \partial_t \tilde{M} \ &\text{in } L^2\left((0,T) \times \Omega_0 \right), \label{68} \\
\eta_{\Delta t} \rightarrow \eta \ \text{in } C\left([0,T];C^1\left(\overline{\Omega_0} \right) \right), \ \ \quad \quad \quad \quad \quad \overline{\eta}_{\Delta t}, \overline{\eta}'_{\Delta t} \rightarrow \eta \ &\text{in } L^\infty \left(0,T;C^1\left(\overline{\Omega_0} \right) \right),\label{2984} \\
\tilde{M}_{\Delta t} \rightarrow \tilde{M} \ \text{in } C\left([0,T];L^p\left(\Omega_0\right) \right),%\quad
\ \ \quad \quad \quad \quad \quad \quad \overline{\tilde{M}}_{\Delta t} \rightarrow \tilde{M} \ &\text{in } L^\infty\left(0,T;L^p\left(\Omega_0\right) \right), \label{3049} \\
\left(\nabla \eta_{\Delta t}\right)^{-1} \rightarrow \left( \nabla \eta\right) \ &\text{in } C\left([0,T] \times \overline{\Omega_0} \right),\label{3w9tiw93ti9w} \\
\left(\nabla \overline{\eta}_{\Delta t}\right)^{-1}, \left(\nabla \overline{\eta}'_{\Delta t}\right)^{-1} \rightarrow \left( \nabla \eta\right)^{-1} \ &\text{in } L^\infty\left(0,T;C \left( \overline{\Omega_0}\right) \right) \label{3036} 
\end{align}

for all $1 \leq p < 6$ and
\begin{align}
\frac{1}{\det \left( \nabla_X \overline{\eta}_{\Delta t} \right)}\overline{\tilde{M}}_{\Delta t} \buildrel\ast\over\rightharpoonup \frac{1}{\det \left( \nabla_X \eta \right)}\tilde{M} \quad \quad \text{in } L^\infty \left(0,T;H^1\left(\Omega_0\right) \right). \label{3042}
\end{align}

Finally, the injectivity of $\overline{\eta}_{\Delta t}(t) \in \mathcal{E}$ in $\Omega_0$ implies the existence of inverse functions $\overline{\eta}_{\Delta t}^{-1}(t,\cdot): \overline{\eta}_{\Delta t}(t,\Omega_0) \rightarrow \Omega_0$ of $\overline{\eta}_{\Delta t}(t,\cdot)$ for all $t \in [0,T]$, for which the convergences \eqref{2984} and \eqref{3036} show that
\begin{align}
\overline{\eta}^{-1}_{\Delta t}(t,\cdot) \rightarrow \eta^{-1} (t,\cdot) \quad \quad \text{pointwise in } \Omega \left(t\right) := \eta \left( t, \Omega_0 \right)\ \text{for all } t \in \left[0,T \right], \label{3054}
\end{align}

where $\eta^{-1}(t,\cdot):\Omega(t) \rightarrow \Omega_0$ denotes the inverse of $\eta(t,\cdot) \in \mathcal{E}$. We point out that with the current configuration $\Omega (t) = \eta (t, \Omega_0)$ of the limit system given, we may also express the magnetization of the limit system in the current configuration,
\begin{align}
M:= M_\eta \left[ \tilde{M} \right] = \frac{1}{\det \left( \left[ \nabla \eta \right] \left( \eta^{-1} \right) \right)} \tilde{M} \left( \eta^{-1} \right). \label{4020}
\end{align}

While the above uniform bounds and convergences hold true on the interval $[0,T]$ for $T>0$ given by Lemma \ref{uniformbounds}, the discrete equation of motion \eqref{3062} is satisfied on this interval only if $T>0$ is in addition chosen such that $\eta_{\Delta t}^k \in \operatorname{int}( \mathcal{E})$ for all $k=1,...,\frac{T}{\Delta t}$. The existence of a time $T>0$ independent of $\Delta t>0$ for which this is indeed true can be seen as another consequence of the uniform bounds \eqref{uniformestimates}, \eqref{3035} and is proved in the following lemma:
\begin{lemma}
\label{injectivity}
There exists a time $T>0$, independent of $\Delta t>0$, such that $\eta_{\Delta t}^k \in \operatorname{int}( \mathcal{E})$ for all $k=1,...,\frac{T}{\Delta t}$. In particular, the equation of motion \eqref{3062} holds true on the interval $[0,T]$ for all test functions $\chi \in \mathcal{D}((0,T)\times \Omega_0)$.
\end{lemma}

\textbf{Proof}

Due to the uniform bounds \eqref{uniformestimates}, \eqref{3035} of $\eta_{\Delta t}$ and the Morrey embedding $H^{1}(0,T)\subset C^{0,\frac{1}{2}}([0,T])$ \cite[Section 1.3.5.8]{novotnystraskraba}), it immediately follows that
\begin{align}
\left\| \eta_{\Delta t} \right\|_{C^{0,\frac{1}{2}}([0,T];H^1(\Omega_0))} \leq c \nonumber
\end{align}

for a constant $c>0$ independent of $\Delta t$. Consequently, it holds that
\begin{align}
\left\| \eta_{\Delta t}\left(t_1\right) - \eta_{\Delta t}\left(t_2\right) \right\|_{H^1(\Omega_0)} \leq c \sqrt{t_1 - t_2} \quad \quad \forall t_1, t_2 \in \left[0,T\right],\ t_1 > t_2. \label{98}
\end{align}

Now let $t_1, t_2 \geq 0$ be such that $t_1 > t_2 + \Delta t$. Then there exist $k,l \in \mathbb{N}$, $k \geq l+1$, such that $t_1 \in ((k-1)\Delta t, k\Delta t]$ and $t_2 \in ((l-1)\Delta t, l\Delta t]$. Consequently,
\begin{align}
(k-l)\Delta t \leq 2\max \left\{ \Delta t, (k-1-l)\Delta t \right\} \leq 2 \left( t_1 - t_2 \right). \label{2988}
\end{align}

Moreover, by definition of the piecewise affine interpolant $\eta_{\Delta t}$ in \eqref{2931} and the piecewise constant interpolant $\overline{\eta}_{\Delta t}$ in \eqref{2930}, we know that $\overline{\eta}_{\Delta t} \left(t_1 \right) = \eta_{\Delta t}\left(k\Delta t \right)$ and $\overline{\eta}_{\Delta t} \left(t_2 \right) = \eta_{\Delta t}\left(l\Delta t \right)$. Therefore, the Hölder continuity \eqref{98} of $\eta_{\Delta t}$ and the estimate \eqref{2988} imply that
\begin{align}
\left\| \overline{\eta}_{\Delta t} \left(t_1 \right) - \overline{\eta}_{\Delta t} \left(t_2 \right) \right\|_{H^1(\Omega_0)} = \left\| \eta_{\Delta t}\left(k\Delta t \right) - \eta_{\Delta t}\left(l\Delta t \right) \right\|_{H^1(\Omega_0)} \leq c \sqrt{(k-l)\Delta t} \leq \sqrt{2}c \sqrt{t_1 - t_2} \label{99}
\end{align}

for all $t_1, t_2 \geq 0$ satisfying $t_1 > t_2 + \Delta t$. Now let $\Gamma>0$ be the constant given by Lemma \ref{injectivityonboundary} in the appendix. From the estimate \eqref{99} we infer the existence of some sufficiently small time $T>0$, independent of $\Delta t$, such that
\begin{align}
\left\| \eta_{\Delta t}^k - \eta_0 \right\|_{H^1(\Omega_0)} < \Gamma \quad \quad \text{for all sufficiently small } \Delta t > 0\text{ and all } k=1,...,\frac{T}{\Delta t}. \nonumber
\end{align}

Hence, Lemma \ref{injectivityonboundary} implies that for all such $k$ the deformation $\eta_{\Delta t}^k$ is injective on $\partial \Omega_0$ and in particular on $N$. 
\begin{comment}
Since $\eta_0 \in \operatorname{int}(\mathcal{E})$ we further know from Remark \ref{interiorpoints} that $\eta_0(N) \bigcap \partial K = \emptyset$. Thus, due to the convergence \eqref{3036} we can also choose $T$ sufficiently small such that
\begin{align}
\eta_{\Delta t}^k (N) \bigcap \partial K = \emptyset \quad \quad \text{for all } k=1,...,\frac{T}{\Delta t}. \nonumber
\end{align}
\end{comment}
From Remark \ref{interiorpoints} it follows that $\eta_{\Delta t}^k \in \operatorname{int}( \mathcal{E})$.
$\hfill \Box$\\

In the following we will first finish our existence proof on the interval $[0,T]$ where $T>0$ with $\frac{T}{\Delta t} \in \mathbb{N}$ is chosen according to Lemma \ref{uniformbounds} and Lemma \ref{injectivity}. Subsequently, we will extend the resulting weak solution to the interval $[0,T')$ where $T'>0$ is chosen as in Theorem \ref{mainresultmagnetoelastic}.\\

Our next goal is the proof of $L^2$-convergence of the stray field, which is needed for the limit passage in both the equation of motion and the magnetic force balance: 

\begin{lemma}
For a suitably chosen subsequence it holds that
\begin{align}
H\left[ \overline{\tilde{M}}_{\Delta t}, \overline{\eta}_{\Delta t} \right] \left( \overline{\eta}_{\Delta t} \right) \rightarrow H\left[ \tilde{M}, \eta \right]\left( \eta \right) \quad \quad \text{in } L^2\left((0,T)\times \Omega_0 \right), \label{3045}
\end{align}

where $H[\tilde{M},\eta] = - \nabla_x \phi[\tilde{M}, \eta]$ and $\phi[\tilde{M}, \eta]$ denotes the solution to the Poisson problem \eqref{3008} associated to $\tilde{M}$ and $\eta$.
\end{lemma}

\textbf{Proof}

Thanks to the uniform bound \eqref{uniformestimates} of $\frac{\overline{\tilde{M}}_{\Delta t}}{\det (\nabla \overline{\eta}_{\Delta t})}$ and the estimate \eqref{3076} for the solution to the Poisson equation \eqref{3008}, we know that
\begin{align}
\left\| H\left[ \overline{\tilde{M}}_{\Delta t}, \overline{\eta}_{\Delta t} \right] \right\|_{L^\infty(0,T;L^2(\mathbb{R}^3))} = \left\| \phi\left[\overline{\tilde{M}}_{\Delta t}, \overline{\eta}_{\Delta t} \right] \right\|_{L^\infty(0,T;\dot{H}^1(\mathbb{R}^3))} \leq c \label{53}
\end{align}

for a constant $c>0$ independent of $\Delta t$. We may thus extract a subsequence such that
\begin{align}
H\left[ \overline{\tilde{M}}_{\Delta t}, \overline{\eta}_{\Delta t} \right] = - \nabla \phi \left[\overline{\tilde{M}}_{\Delta t}, \overline{\eta}_{\Delta t} \right] \rightharpoonup - \nabla \phi \left[ \tilde{M}, \eta \right] = H \left[ \tilde{M}, \eta \right] \quad \quad \text{in } L^2\left(\left(0,T \right) \times \mathbb{R}^3 \right),\label{50}
\end{align}

where the limit function could be identified as the (negative) gradient of the solution $\phi [\tilde{M}, \eta](t,\cdot) \in \dot{H}^1(\mathbb{R}^3)$ to the Poisson equation \eqref{3008} associated to $M$ and $\eta$ by passing to the limit in the corresponding Poisson equation satisfied by $\phi [\overline{\tilde{M}}_{\Delta t}, \overline{\eta}_{\Delta t}]$ under exploitation of the uniform convergence \eqref{2984} of the deformation.

In order to strengthen \eqref{50} to strong convergence, we use the definition of the stray field via \eqref{3008}, and estimate
\begin{align}
&\left| \int_0^{T} \int_{\mathbb{R}^3} \left| H\left[ \overline{\tilde{M}}_{\Delta t}, \overline{\eta}_{\Delta t} \right] \right|^2 \ dxdt - \int_0^{T} \int_{\mathbb{R}^3} \left| H\left[ \tilde{M}, \eta \right] \right|^2 \ dxdt \right| \nonumber \\
\leq & \left| \int_0^{T} \int_{\Omega_0} \left( \overline{\tilde{M}}_{\Delta t} - \tilde{M} \right) \cdot H\left[ \overline{\tilde{M}}_{\Delta t}, \overline{\eta}_{\Delta t} \right] \left( \overline{\eta}_{\Delta t} \right)\ dXdt \right| \nonumber \\
&+ \left| \int_0^{T} \int_{\Omega_0} \tilde{M} \cdot \left( H\left[ \overline{\tilde{M}}_{\Delta t}, \overline{\eta}_{\Delta t} \right] \left( \overline{\eta}_{\Delta t} \right) - H\left[ \tilde{M}, \eta \right]\left( \overline{\eta}_{\Delta t} \right) \right) \ dXdt \right| \nonumber \\
&+ \left| \int_0^{T} \int_{\Omega_0} \tilde{M} \cdot \left( H\left[ \tilde{M}, \eta \right] \left( \overline{\eta}_{\Delta t} \right) - H\left[ \tilde{M}, \eta \right]\left( \eta \right) \right) \ dXdt \right|. \label{51}
\end{align}

Here, the first integral on the right-hand side tends to zero due to the strong convergence \eqref{3049} of the magnetization, while the third integral vanishes thanks to the uniform convergence \eqref{2984} of the deformation and the uniform bound \eqref{uniformestimates} of its Jacobian away from zero. In order to show that also the second integral on the right-hand side of the inequality \eqref{51} vanishes we choose a sequence of functions $(\tilde{M}_n)_{n \in \mathbb{N}} \subset \mathcal{D}((0,T)\times \Omega_0)$ such that $\tilde{M}_n \rightarrow \tilde{M}$ in $L^6((0,T)\times \Omega_0)$. Then we estimate, via the triangle inequality,
\begin{align}
&\left| \int_0^{T} \int_{\Omega_0} \left(\tilde{M} - \tilde{M}_n + \tilde{M}_n \right) \cdot \left( H\left[ \overline{\tilde{M}}_{\Delta t}, \overline{\eta}_{\Delta t} \right] \left( \overline{\eta}_{\Delta t} \right) - H\left[ \tilde{M}, \eta \right]\left( \overline{\eta}_{\Delta t} \right) \right) \ dXdt \right| \nonumber \\
\leq& c\left\| \tilde{M} - \tilde{M}_n \right\|_{L^2((0,T)\times \Omega_0)} + \left| \int_0^{T} \int_{\overline{\eta}_{\Delta t}(t,\Omega_0)} \frac{\tilde{M}_n \left(\left( \overline{\eta}_{\Delta t} \right)^{-1} \right)}{\det\left(\left[\nabla \overline{\eta}_{\Delta t} \right] \left( \overline{\eta}_{\Delta t}^{-1} \right) \right)} \cdot \left( H\left[ \overline{\tilde{M}}_{\Delta t}, \overline{\eta}_{\Delta t} \right] - H\left[ \tilde{M}, \eta \right] \right) \ dxdt \right|. \nonumber
\end{align}

Making use of the pointwise convergence \eqref{3054} of the inverse deformation, the Vitali convergence theorem and the weak convergence \eqref{50} of the stray field, we see that the right-hand side of this inequality vanishes for $\Delta t \rightarrow 0$ and $n \rightarrow \infty$. Consequently, the right-hand side of the estimate \eqref{51} also tends to zero for $\Delta t \rightarrow 0$. This shows convergence of the $L^2((0,T)\times \mathbb{R}^3)$-norm of $H[\overline{\tilde{M}}_{\Delta t}, \overline{\eta}_{\Delta t}]$ to the one of $H[\tilde{M}, \eta]$. In combination with the weak convergence \eqref{50} of $H[\overline{\tilde{M}}_{\Delta t}, \overline{\eta}_{\Delta t}]$, this implies strong convergence of $H[\overline{\tilde{M}}_{\Delta t}, \overline{\eta}_{\Delta t}]$ in $L^2((0,T)\times \mathbb{R}^3)$ and hence the desired strong convergence \eqref{3045}.
$\hfill \Box$\\

We can now turn to the limit passage in the discrete magnetic force balance \eqref{3034}. For the anisotropy energy density, onto which we imposed the continuity and boundedness assumptions \eqref{3063}, \eqref{3064}, we infer from the strong convergences \eqref{2984}, \eqref{3049} and the Vitali convergence theorem that
\begin{align}
\tilde{\Psi}_M \left( \nabla \overline{\eta}_{\Delta t}, \overline{\tilde{M}}_{\Delta t} \right) \rightarrow \tilde{\Psi}_M \left( \nabla \eta, \tilde{M} \right)\quad \quad \text{in } L^p \left((0,T) \times 
\Omega_0 \right) \quad \text{for some } p>\frac{6}{5}. \label{2928}
\end{align}

Then, we test the discrete magnetic force balance \eqref{3034} by functions of the form $\det( \nabla \overline{\eta}_{\Delta t}) \tilde{\hat{M}} \in L^\infty(0,T;H^1(\Omega_0))$ with $\tilde{\hat{M}} \in L^\infty ( 0,T;H^1(\Omega_0))$. In the resulting identity we pass to the limit by combining the convergence \eqref{2928} with the uniform convergences \eqref{2984}, \eqref{3036} of $\overline{\eta}_{\Delta t}$, $\overline{\eta}'_{\Delta t}$ and $(\nabla \overline{\eta}_{\Delta t})^{-1}$, the weak and strong convergences \eqref{68} and \eqref{3049} of $\overline{\tilde{M}}_{\Delta t}$ and $\tilde{M}_{\Delta t}$ and the convergence \eqref{3045} of $H[\overline{\tilde{M}}_{\Delta t}, \overline{\eta}_{\Delta t}](\overline{\eta}_{\Delta t})$. This yields the limit equation
\begin{align}
& \int_0^T \left\langle D_2\tilde{E} \left(\eta, \tilde{M} \right), \tilde{\hat{M}} \right\rangle_{\left( H^1 \left(\Omega_0 \right) \right)^* \times H^1 \left(\Omega_0 \right)} + \left\langle D_3\tilde{R} \left( \eta, \partial_t \eta, \partial_t \tilde{M} \right), \tilde{\hat{M}} \right\rangle_{L^2(\Omega_0) \times L^2(\Omega_0)} \ dt = 0 \label{2981} 
\end{align}

for all $\tilde{\hat{M}} \in L^\infty(0,T;H^1(\Omega_0))$. We test this equation by $\hat{M}(\eta) \in L^\infty(0,T;H^1(\Omega_0))$ for some arbitrary test function $\hat{M} \in L^\infty (0,T;H^1(\Omega (\cdot)))$. Transforming the resulting identity to the current configuration we infer that the pair $(\eta, M)$, where $M$ denotes the magnetization in the current configuration defined in \eqref{4020}, satisfies the desired magnetic force balance \eqref{3057}.\\

It remains to pass to the limit in the discrete equation of motion \eqref{3062}. In order to pass to the limit in the exchange energy terms we first show strong convergence of the magnetization gradient. To this end we make use of Minty's trick: We return to the discrete magnetic force balance \eqref{3034} and test it by the test function $\tilde{\hat{M}} = \overline{\tilde{M}}_{\Delta t} \in L^\infty (0,T;H^1(\Omega_0))$. In the resulting identity we pass to the limit once more (cf.\ the derivation of \eqref{2981}) and compare the result to the magnetic force balance \eqref{2981} tested by $\frac{1}{\det(\nabla \eta)}\tilde{M}\in L^\infty (0,T;H^1(\Omega))$. This leads to the relation
\begin{align}
&\lim_{\Delta t \rightarrow 0} \int_0^{T}\int_{\Omega_0} 2A\det \left( \nabla \overline{\eta}_{\Delta t} \right) \left| \nabla \left( \frac{1}{\det \left( \nabla \overline{\eta}_{\Delta t} \right)}\overline{\tilde{M}}_{\Delta t} \right) \left( \nabla \overline{\eta}_{\Delta t} \right)^{-1} \right|^2 \ dXdt \nonumber \\
=& \int_0^{T}\int_{\Omega_0} 2A \det \left( \nabla \eta \right) \left| \nabla \left( \frac{1}{\det \left( \nabla \eta \right)}\tilde{M} \right) \left( \nabla \eta \right)^{-1} \right|^2\ dXdt. \nonumber
\end{align}

In combination with the weak convergence \eqref{3042}, the bound \eqref{uniformestimates} of $\det (\nabla \overline{\eta}_{\Delta t})$ away from zero and the uniform convergence \eqref{2984} of both $\overline{\eta}_{\Delta t}$ and its gradient this implies the desired strong convergence
\begin{align}
\nabla \left( \frac{1}{\det \left( \nabla \overline{\eta}_{\Delta t} \right)}\overline{\tilde{M}}_{\Delta t} \right) \rightarrow \nabla \left( \frac{1}{\det \left( \nabla \eta \right)}\tilde{M} \right) \quad \quad \text{in } L^2\left(\left(0,T\right)\times \Omega_0 \right). \label{3058}
\end{align}

In order to pass to the limit in the exchange energy terms (stemming from the micromagnetic energy \eqref{4041}) in the equation of motion \eqref{3062}, we combine the convergence \eqref{3058} with the weak convergence \eqref{40} of $\nabla^2 \overline{\eta}_{\Delta t}$ and the strong convergence \eqref{3049} of $\overline{\tilde{M}}_{\Delta t}$, to deduce the convergence
\begin{align}
\nabla \left( \frac{\operatorname{tr} \left( \nabla \chi \left( \nabla \overline{\eta}_{\Delta t} \right)^{-1} \right)}{\det \left( \nabla \overline{\eta}_{\Delta t} \right)} \overline{\tilde{M}}_{\Delta t} \right) \rightharpoonup \nabla \left( \frac{\operatorname{tr} \left( \nabla \chi \left( \nabla \eta \right)^{-1} \right)}{\det \left( \nabla \eta \right)} \tilde{M} \right) \quad \quad \text{in } L^2 \left(\left(0,T \right)\times \Omega_0 \right) \label{exenconv}
\end{align}

for any $\chi \in \mathcal{D}((0,T)\times \Omega_0)$. The latter two convergences allow us to indeed pass to the limit in the exchange energy terms in the equation of motion \eqref{3062},
\begin{align}
&\int_0^T \int_{\Omega_0} A \det\left( \nabla \overline{\eta}_{\Delta t} \right) \left| \nabla \left( \frac{1}{\det \left( \nabla \overline{\eta}_{\Delta t} \right)} \overline{\tilde{M}}_{\Delta t} \right) \left( \nabla \overline{\eta}_{\Delta t} \right)^{-1} \right|^2 \left(\left( \nabla \overline{\eta}_{\Delta t} \right)^{-1}\right)^T : \nabla \chi \nonumber \\
&-2A \det\left( \nabla \overline{\eta}_{\Delta t} \right) \left[ \nabla \left( \frac{1}{\det \left( \nabla \overline{\eta}_{\Delta t} \right)} \overline{\tilde{M}}_{\Delta t} \right) \left( \nabla \overline{\eta}_{\Delta t} \right)^{-1} \right] : \left[ \nabla \left( \frac{\operatorname{tr} \left( \nabla \chi \left( \nabla \overline{\eta}_{\Delta t} \right)^{-1} \right)}{\det \left( \nabla \overline{\eta}_{\Delta t} \right)} \overline{\tilde{M}}_{\Delta t} \right) \right. \left( \nabla \overline{\eta}_{\Delta t} \right)^{-1} \nonumber \\
&\left. \nabla \left( \frac{1}{\det \left( \nabla \overline{\eta}_{\Delta t} \right)} \overline{\tilde{M}}_{\Delta t} \right) \left( \nabla \overline{\eta}_{\Delta t} \right)^{-1} \nabla \chi \left( \nabla \overline{\eta}_{\Delta t} \right)^{-1} \vphantom{\left( \frac{\operatorname{tr} \left( \nabla \chi \left( \nabla \overline{\eta}_{\Delta t} \right)^{-1} \right)}{\det \left( \nabla \overline{\eta}_{\Delta t} \right)} \overline{\tilde{M}}_{\Delta t} \right)} \right] \ dXdt \nonumber \\
\rightarrow & \int_0^T \int_{\Omega_0} A \det\left( \nabla \eta \right) \left| \nabla \left( \frac{1}{\det \left( \nabla \eta \right)} \tilde{M} \right) \left( \nabla \eta \right)^{-1} \right|^2 \left(\left( \nabla \eta \right)^{-1}\right)^T : \nabla \chi \nonumber \\
&-2A \det\left( \nabla \eta \right) \left[ \nabla \left( \frac{1}{\det \left( \nabla \eta \right)} \tilde{M} \right) \left( \nabla \eta \right)^{-1} \right] : \left[ \nabla \left( \frac{\operatorname{tr} \left( \nabla \chi \left( \nabla \eta \right)^{-1} \right)}{\det \left( \nabla \eta \right)} \tilde{M} \right) \left( \nabla \eta \right)^{-1} \right. \nonumber \\
&+ \left. \vphantom{\nabla \left( \frac{\operatorname{tr} \left( \nabla \chi \left( \nabla \eta \right)^{-1} \right)}{\det \left( \nabla \eta \right)} \tilde{M} \right) \left( \nabla \eta \right)^{-1}} \nabla \left( \frac{1}{\det \left( \nabla \eta \right)} \tilde{M} \right) \left( \nabla \eta \right)^{-1} \nabla \chi \left( \nabla \eta \right)^{-1} \right] \ dXdt. \label{exchangelimit}
\end{align}

For the limit passage in the stray field term we require (local) weak convergence of the gradient of the stray field. In order to show this we consider an arbitrary compact set $K \subset \Omega_0$. From the compactness of $K$ and the uniform convergence \eqref{2984} of $\overline{\eta}_{\Delta t}$ we infer the existence of some constant $\delta >0$ independent of $\Delta t$ such that
\begin{align}
\text{dist} \left( \overline{\eta}_{\Delta t}\left(t,K\right), \partial \left( \overline{\eta}_{\Delta t}(t,\Omega_0) \right) \right) \geq \delta \quad \quad \forall t \in [0,T] \nonumber
\end{align}

for all sufficiently small $\Delta t > 0$. This allows us to make use of the bound \eqref{78}, given by Lemma \ref{poissonequation} in the appendix for solutions to the Poisson equation, and infer that
\begin{align}
&\left\| \nabla H \left[ \overline{\tilde{M}}_{\Delta t}, \overline{\eta}_{\Delta t} \right] \left(\overline{\eta}_{\Delta t} \right) \right\|_{L^2((0,T)\times K)}^2 \nonumber \\
\leq& c \left\| \nabla_x H \left[ \overline{\tilde{M}}_{\Delta t}, \overline{\eta}_{\Delta t} \right] \right\|^2_{L^2 (0,T;L^2(\overline{\eta}_{\Delta t} (\cdot, K)))} \leq c \int_0^{T}\int_{\overline{\eta}_{\Delta t}(t,\Omega_0)} \left| \overline{M}_{\Delta t} \right|^2 + \left| \nabla_x \overline{M}_{\Delta t} \right|^2 \ dxdt \leq c \nonumber
\end{align}

for a constant $c>0$ independent of $\Delta t$ thanks to the uniform bounds \eqref{uniformestimates}. Thus, using a diagonal argument, we may extract another subsequence for which it holds that
\begin{align}
\nabla H \left[ \overline{\tilde{M}}_{\Delta t}, \overline{\eta}_{\Delta t} \right] \left(\overline{\eta}_{\Delta t} \right) \rightharpoonup \nabla H \left[ \tilde{M}, \eta \right] \left(\eta \right) \quad \quad \text{in } L^2\left(\left(0,T\right) \times K \right) \label{3070}
\end{align}

for any compact set $K \subset \Omega_0$, where the identification of the limit function results from the already known convergence \eqref{3045}. Since the test functions $\chi \in \mathcal{D}((0,T)\times \Omega_0)$ in the momentum equation are compactly supported, the convergence \eqref{3070} suffices to pass to the limit in the stray field term (stemming from the micromagnetic energy \eqref{4041}) in the approximate equation of motion \eqref{3062},
\begin{align}
&\int_0^T \int_{\Omega_0} \mu \left[ \left( \left( \nabla H\left[ \overline{\tilde{M}}_{\Delta t}, \overline{\eta}_{\Delta t} \right] \left( \overline{\eta}_{\Delta t} \right) \right) \left(\nabla \overline{\eta}_{\Delta t} \right)^{-1} \right)^T \overline{\tilde{M}}_{\Delta t} \right] \cdot \chi \ dxdt \nonumber \\
\rightarrow & \int_0^T \int_{\Omega_0} \mu \left[ \left( \left( \nabla H\left[ \tilde{M}, \eta \right] \left( \eta \right) \right) \left(\nabla \eta \right)^{-1} \right)^T \tilde{M} \right] \cdot \chi \ dxdt. \label{strayfieldlimit}
\end{align}

For the anisotropy energy density $\tilde{\Psi}$ and the elastic energy density $W$, we recall the continuity and boundedness assumptions \eqref{3063}, \eqref{3090}, \eqref{2980}, which, in combination with the strong convergences \eqref{2984} and \eqref{3049} and the Vitali convergence theorem yield
\begin{align}
W' \left( \nabla \overline{\eta}_{\Delta t} \right) \rightarrow W' \left( \nabla \eta \right) \quad \quad &\text{in } L^p\left((0,T)\times \Omega_0 \right), \label{3065} \\
\tilde{\Psi}_F \left( \nabla \overline{\eta}_{\Delta t}, \overline{\tilde{M}}_{\Delta t} \right) \rightarrow \tilde{\Psi}_F \left( \nabla \eta, \tilde{M} \right)\quad \quad &\text{in } L^p\left((0,T)\times \Omega_0 \right) \label{3065.5}
\end{align}

for some $p>1$. Finally, we remark that due to the uniform bound \eqref{uniformestimates} we can extract another subsequence and find some limit function $z \in L^\infty (0,T;L^\frac{q}{q-1}(\Omega_0))$ such that
\begin{align}
\left| \nabla^2 \overline{\eta}_{\Delta t} \right|^{q-2}\nabla^2 \overline{\eta}_{\Delta t} \buildrel\ast\over\rightharpoonup z \quad \quad \text{in } L^\infty \left(0,T;L^\frac{q}{q-1}\left(\Omega_0\right) \right). \label{3072}
\end{align}

In order to identify the limit function $z$, we test (after a density argument) the approximate equation of motion \eqref{3062} by $(\overline{\eta}_{\Delta t} - \eta) \phi$ for an arbitrary non-negative smooth cutoff function $\phi \in \mathcal{D}((0,T)\times \Omega_0)$. Combining the convergences \eqref{3058}--\eqref{3072}, we pass to the limit in the resulting identity and deduce that
\begin{align}
\int_0^{T} \int_{\Omega_0} \left( \left|\nabla^2 \overline{\eta}_{\Delta t} \right|^{q-2}\nabla^2 \overline{\eta}_{\Delta t} \right)\ \vdots \ \nabla^2 \left( \left(\overline{\eta}_{\Delta t} - \eta \right) \phi \right) \ dXdt \rightarrow 0, \label{convergence}
\end{align}

where the triple dot product $\vdots$ denotes the generalization of the classical dot product to an inner product between $3\times 3 \times 3$-tensors. The convergence \eqref{convergence} in combination with the weak convergence \eqref{40} of $\overline{\eta}_{\Delta t}$ shows that %cf.\ [variational approach to fsi, page 28] or iphone photos 18.08.2023
\begin{align}
&\left\| \left( \nabla^2 \overline{\eta}_{\Delta t} - \nabla_X^2 \eta \right) \phi^\frac{1}{q} \right\|_{L^q((0,T)\times \Omega_0)}^q \nonumber \\
\leq& c \int_0^{T} \int_{\Omega_0} \left( \left|\nabla^2 \overline{\eta}_{\Delta t} \right|^{q-2}\nabla^2 \overline{\eta}_{\Delta t} - \left|\nabla^2 \eta \right|^{q-2}\nabla^2 \eta \right)\ \vdots \ \nabla^2 \left(\overline{\eta}_{\Delta t} - \eta \right) \phi \ dXdt \rightarrow 0. \nonumber
\end{align}

Since $\phi$ can be chosen as an approximation of the constant function $1$ on $(0,T)\times \Omega_0$, this relation is sufficient to identify the limit function $z$ in the convergence \eqref{3072} as
\begin{align}
z = \left|\nabla^2 \eta \right|^{q-2}\nabla^2 \eta \quad \quad \text{a.e. in } \left(0,T\right) \times \Omega_0. \label{3074}
\end{align}

Now, combining the convergences \eqref{exchangelimit}, \eqref{strayfieldlimit}--\eqref{3072} and the identity \eqref{3074} we may pass to the limit in the approximate equation of motion \eqref{3062} for arbitrary test functions $\chi \in \mathcal{D}((0,T)\times \Omega_0)$ and infer the desired equation of motion \eqref{3075}.\\

Having shown the equation of motion \eqref{3075} and the magnetic force balance \eqref{3057} (cf.\ the identity \eqref{2981}), we have so far proved the existence of a weak solution in the sense of Definition \ref{weaksolutionsmagnetoelastic} on the interval $[0,T)$ with $T>0$ chosen according to Lemma \ref{uniformbounds} and Lemma \ref{injectivity}. Indeed, the regularity \eqref{2958}, \eqref{2957} of the deformation and the magnetization follows from the convergences \eqref{40}, \eqref{68} and the fact that the set $\mathcal{E}$ is closed with respect to weak convergence in $W^{2,q}(\Omega_0)$. From the uniform in time convergences \eqref{2984} and \eqref{3049} of the deformation and the magnetization we further conclude the initial conditions \eqref{3100}, \eqref{3101}. The $L^\infty (0,T;H^1(\Omega (\cdot)))$-regularity \eqref{2951} of $M$ follows from the $L^\infty (0,T;H^1(\Omega_0))$-regularity of $\frac{1}{\det(\nabla \eta)}\tilde{M}$ (cf.\ \eqref{3042}), the essential bounds of $\det (\nabla \eta)$ and $(\nabla \eta)^{-1}$ (cf.\ \eqref{2984}, \eqref{3036}) and a transformation to the reference configuration. \\

Finally, it remains to show that the existence time $T$ can be chosen as $T=T'$ for $T'>0$ as in Theorem~\ref{mainresultmagnetoelastic}. To this end we denote by $T_{\text{max}}> 0$ the maximal time such that $(\eta, \tilde{M})$ is a weak solution on each interval $[0,T)$ with $0 < T < T_{\text{max}}$. We assume that $T_{\text{max}} < \infty$ and $\liminf_{t \rightarrow T_{\operatorname{max}}}\tilde{E}(\eta (t),\tilde{M}(t))< \infty.$ Under exploitation of the latter bound, the functions $(\eta, \tilde{M})$ can be extended to the time $T_{\text{max}}$ by values $\eta (T_{\operatorname{max}}) \in \mathcal{E}$ and $\tilde{M}(T_{\operatorname{max}})\in H^1(\Omega_0)$ such that
\begin{align}
\eta &\in L^\infty \left(0,T_{\operatorname{max}};\mathcal{E} \right) \bigcap C\left(\left[0,T_{\operatorname{max}}\right];C^1\left(\overline{\Omega_0} \right) \right), \label{2939} \\
\tilde{M} &\in L^\infty \left(0,T_{\operatorname{max}};H^1\left(\Omega_0 \right) \right) \bigcap C\left(\left[0,T_{\operatorname{max}}\right];L^2\left(\Omega_0 \right) \right), \label{2940} \\
\partial_t \eta &\in L^2\left(0, T_{\text{max}};H^1\left(\Omega_0\right)\right),\quad \quad \partial_t \tilde{M} \in L^2\left(\left(0, T_{\text{max}}\right)\times \Omega_0 \right). \label{2941}
\end{align}

With the regularity \eqref{2939}--\eqref{2941} at hand we may conclude the proof via a contradiction argument: We assume in addition that there occurs no self-contact of the material at the time $T_{\text{max}}$, i.e. $\eta(T_{\text{max}}) \in \text{int}( \mathcal{E})$ (cf.\ Remark \ref{selfcontact}). Then we may apply the local existence result we have proved so far to construct a solution on the interval $[T_{\text{max}}, T_{\text{max}} + \epsilon)$ for some small $\epsilon > 0$. Due to the regularity \eqref{2939}--\eqref{2941} the two solutions can be assembled to a solution on the interval $[0,T_{\text{max}} + \epsilon)$. This, however, results in a contradiction to the maximality of $T_{\text{max}}$. Therefore, the assumption $\eta(T_{\text{max}}) \in \text{int}( \mathcal{E})$ is wrong. It follows that there occurs a self-contact of the material at $T_{\text{max}}$ and hence $T_{\text{max}}$ coincides with the ending time $T'$ as specified in Theorem \ref{mainresultmagnetoelastic}. This finishes the proof of Theorem \ref{mainresultmagnetoelastic}.

\appendix

\section{Appendix}

In the appendix we summarize several auxiliary results necessary for the proof of our main result Theorem \ref{mainresultmagnetoelastic}. For the existence and the uniform bounds of the stray field we make use of the following classical results for the Poisson equation, cf.\ \cite[Theorem 5.1]{praetorius}, \cite[Theorem 8.8]{gilbarg}.

\begin{lemma} \label{poissonequation}

Let $\Omega \subset \mathbb{R}^3$ be a bounded domain. Let $M \in H^1(\Omega)$ be extended by $0$ outside of $\Omega$. Then there exists a unique solution $\phi \in \dot{H}^1(\mathbb{R}^3)$ to the Poisson equation
\begin{align}
\int_{\mathbb{R}^3} \nabla_x \phi \cdot \nabla_x \psi \ dx = \int_{\Omega} M \cdot \nabla_x \psi \ dx \quad \quad \forall \psi \in \dot{H}^1\left(\mathbb{R}^3 \right). \label{77}
\end{align}

Further there exists a constant $c>0$, independent of $\Omega$, such that
\begin{align}
\left\| \phi \right\|_{\dot{H}^1(\mathbb{R}^3)} = \left\| \nabla_x \phi \right\|_{L^2(\mathbb{R}^3)} \leq c\|M\|_{L^2(\Omega)} \label{3076}
\end{align}

Moreover, for any $\delta > 0$ there exists a constant $c(\delta)>0$ independent of $M$ and $\Omega$ such that
\begin{align}
\left\| \nabla_x^2 \phi \right\|_{L^2(K)} &\leq c(\delta)\|M\|_{H^1(\Omega)} \label{78}
\end{align}

for all compact subsets $K \subset \Omega$ satisfying $\operatorname{dist}(K, \partial \Omega) \geq \delta$.

\end{lemma}

\textbf{Proof}

The existence of a unique solution $\phi$ to the Poisson equation \eqref{77} as well as its $\dot{H}^1(\mathbb{R}^3)$-bound \eqref{3076} is proved in \cite[Theorem 5.1]{praetorius}. A proof for a local $H^2$-bound in terms of the $H^1(\Omega)$-norm of $M$ and the $L^2(\Omega)$-norm of $\nabla_x \phi$ is given in \cite[Theorem 8.8]{gilbarg}, so the local $L^2$-bound \eqref{78} of $\nabla_x^2 \phi$ follows immediately from the $\dot{H}^1(\mathbb{R}^3)$-bound \eqref{3076}. We point out that, strictly speaking, the bound of the second gradient in \cite[Theorem 8.8]{gilbarg} depends on the $H^1(\Omega)$-norm of $\phi$ instead of only the $L^2(\Omega)$-norm of $\nabla_x \phi$. However, \cite[Theorem 8.8]{gilbarg} is formulated for the case of more general elliptic equations and a look into its proof shows that boundedness of the $L^2(\Omega)$-norm of $\nabla_x \phi$ is indeed sufficient for the setting of the Poisson equation.
$\hfill \Box$\\

We proceed by summarizing miscellaneous auxiliary results concerning our use of a time discretization. We begin with the following variant of \cite[Theorem 8.9]{roubicek}, which is used in the limit passage from the discretized system back to the continuous in time setting, in order to guarantee that the weak limits of different interpolants of the same discrete functions coincide.
\begin{lemma}
\label{equalityofrothelimits}
Let $T>0$, $\Delta t > 0$, with $\frac{T}{\Delta t} \in \mathbb{N}$ and let $\Omega \subset \mathbb{R}^3$ be a domain. Let further $h_{\Delta t}^k:\Omega \rightarrow \mathbb{R}^l$, $k=0,...,\frac{T}{\Delta t}$, $l \in \mathbb{N}$, be time-independent functions with piecewise affine and piecewise constant interpolants
\begin{align}
h_{\Delta t}(t) &:= \left( \frac{t}{\Delta t} - (k-1)\right)h_{\Delta t}^k + \left(k - \frac{t}{\Delta t}\right) h_{\Delta t}^{k-1}\ \ \ &\text{for }& (k-1)\Delta t < t \leq k \Delta t,\quad k = 1,...,\frac{T}{\Delta t}, \nonumber \\
\overline{h}_{\Delta t}(t) &:= h_{\Delta t}^k \ \ \ &\text{for }& (k-1)\Delta t < t \leq k \Delta t,\quad k = 0,...,\frac{T}{\Delta t}, \nonumber \\
\overline{h}'_{\Delta t}(t) &:= h_{\Delta t}^{k-1} \ \ \ &\text{for }& (k-1)\Delta t < t \leq k \Delta t,\quad k = 1,...,\frac{T}{\Delta t}. \nonumber
\end{align}

Assume moreover that
\begin{align}
h_{\Delta t} \buildrel\ast\over\rightharpoonup h \quad \text{in } L^\infty(0,T;L^2(\Omega)),\quad \overline{h}_{\Delta t} \buildrel\ast\over\rightharpoonup \overline{h} \quad \text{in } L^\infty(0,T;L^2(\Omega)),\quad \overline{h}'_{\Delta t} \buildrel\ast\over\rightharpoonup \overline{h}' \quad \text{in } L^\infty(0,T;L^2(\Omega)). \nonumber
\end{align}

Then it holds that
\begin{align}
h = \overline{h} = \overline{h}'. \label{918}
\end{align}

\end{lemma}

\textbf{Proof}

The proof, which is performed by comparing the limit of the functions $\overline{h}_{\Delta t}$, $\overline{h}'_{\Delta t}$ to the one of $h_{\Delta t}$ in the pairing with piecewise constant in time functions respectively, can be found in the proof of \cite[Theorem 8.9]{roubicek}, cf.\ also \cite[Lemma A.3.1]{thesis}.
$\hfill \Box$\\

When passing to the limit in the time discretization we further face the situation of having to pass to the limit in discretized versions of given time-dependent functions. In order to deal with this we use the following version of \cite[Lemma 8.7]{roubicek}.
\begin{lemma}\label{rothelimitextfunctions}\ \\
Let $T>0$, $\Delta t>0$ with $\frac{T}{\Delta t} \in \mathbb{N}$ and let $\Omega \subset \mathbb{R}^3$ be a domain.
\begin{itemize}
\item[(i)] Let $\gamma = \gamma ( \Delta t)>0,$ $\gamma (\Delta t)\rightarrow 0$ for $\Delta t \rightarrow 0$. Let $h \in L^\infty((0,T)\times \Omega)$ and define
\begin{align}
h_\gamma (t):= \int_0^T \theta _\gamma \left( t + \xi_\gamma (t) - s \right)h(s)\ ds,\quad \quad \xi_\gamma (t) := \gamma \frac{T-2t}{T}, \nonumber
\end{align}

where $\theta_\gamma :\mathbb{R}\to \mathbb{R}$ denotes a mollifier with support in $[-\gamma, \gamma]$. Set further
\begin{align}
\overline{h}_{\Delta t}(t) := h_{\gamma (\Delta t)}(k\Delta t)\quad \quad \text{for } (k-1)\Delta t < t \leq k\Delta t,\quad k=1,...,\frac{T}{\Delta t}. \nonumber
\end{align}

Then
\begin{align}
\overline{h}_{\Delta t} \rightarrow h \quad \text{in } L^p((0,T)\times \Omega) \quad \quad \forall 1 \leq p < \infty. \nonumber
\end{align}

\item[(ii)] Let $h \in W^{l,q}([0,T];W^{j,p}(\mathbb{R}^3))$, $l,j \in \mathbb{N}_0$, $1 \leq p,q < \infty$ and define
\begin{align}
\overline{h}_{\Delta t}(t) := \frac{1}{\Delta t} \int_{(k-1)\Delta t}^{k \Delta t} h(s)\ ds \quad \quad \text{for } (k-1)\Delta t < t \leq k\Delta t,\quad k=1,...,\frac{T}{\Delta t}. \nonumber
\end{align}

Then
\begin{align}
\overline{h}_{\Delta t} \rightarrow h \quad \text{in } C\left([0,T];W^{j,p}(\mathbb{R}^3)\right). \nonumber
\end{align}

\end{itemize}
\end{lemma}

\textbf{Proof}

The statement (i) is proved in \cite[Lemma 8.7]{roubicek}; for a proof of the statement (ii), which follows by similar arguments, we refer to \cite[Lemma A.3.2]{thesis}.
$\hfill \Box$\\

In the derivation of the discrete energy estimate in Section \ref{proofsection} we control the discrete difference quotient of the external magnetic field $H_{\text{ext}}$ through its classical time derivative $\partial_t H_{\text{ext}}$ via the following result, which constitutes a version of the estimate (8.72) in the proof of \cite[Theorem 8.18]{roubicek}.
\begin{lemma} \label{differencequotientestimate}
Let $\Delta t > 0$, let $H_{\operatorname{ext}} \in W^{1,\frac{4}{3}}(0,\infty;W^{1,\frac{4}{3}}(\mathbb{R}^3 ))$ and set
\begin{align}
(H_{\operatorname{ext}})_{\Delta t}^k := \frac{1}{\Delta t} \int_{(k-1)\Delta t}^{k\Delta t} H_{\operatorname{ext}} (t)\ dt \quad \quad \forall k \in \mathbb{N}. \nonumber
\end{align}

Then
\begin{align}
\Delta t \sum_{l=2}^k \left\| \frac{\left(H_{\text{ext}} \right)_{\Delta t}^l - \left(H_{\text{ext}} \right)_{\Delta t}^{l-1}}{\Delta t} \right\|_{L^\frac{4}{3}(\mathbb{R}^3)}^\frac{4}{3} \leq \left\| \partial_t H_{\text{ext}} \right\|_{L^\frac{4}{3}((0,\infty)\times \mathbb{R}^3)}^\frac{4}{3} \quad \quad \forall k \in \mathbb{N}. \nonumber
\end{align}
\end{lemma}

\textbf{Proof}

The statement can be proved via an application of Jensen's inequality in the same way as the corresponding estimate (8.72) in the proof of \cite[Theorem 8.18]{roubicek}. The details can be found in \cite[Lemma A.3.4]{thesis}.
$\hfill \Box$\\

We further summarize several auxiliary results regarding the deformability of the domain. For the proof of the main result Theorem \ref{mainresultmagnetoelastic} it is crucial to know that the determinants of the gradients of deformations with uniformly bounded energy are uniformly bounded away from zero. Such a bound is proved in \cite[Section 2.3]{variationalapproachtofsi} (cf.\ also \cite[Lemma A.7.1]{thesis}), the proof therein in turn follows ideas from \cite{healeykromer}. For the convenience of the reader we restate the result in the following lemma.

\begin{lemma} \label{boundawayfromzero}
Let $E_0>0$ be given. Then there exists a constant $c = c(E_0)>0$ such that for all deformations
\begin{align}
\eta \in \mathcal{E} \quad \text{with} \quad \tilde{E}_{\text{el}}(\eta) \leq E_0 \label{2955}
\end{align}

it holds that
\begin{align}
\det \left(\nabla \eta \right) \geq c \quad \text{in } \Omega_0. \nonumber
\end{align}
\end{lemma}

\textbf{Proof}

The proof of this result, which is based upon the presence of the quantities $\frac{1}{( \det ( \nabla \eta ))^a}$ and $\frac{1}{q} \left| \nabla^2 \eta \right|^q$ in the elastic energy, can be found in detail in \cite[Section 2.3]{variationalapproachtofsi} (cf.\ also \cite[Lemma A.7.1]{thesis}).
$\hfill \Box$\\

Moreover, in order to find a small time interval $[0,T]$ on which a deformation $\eta (t)$, $t \in [0,T]$, remains injective on $\partial \Omega_0$ provided that $\eta(0)$ is injective on $\partial \Omega_0$, we use the following version of \cite[Lemma 2.5, Proposition 2.7]{variationalapproachtofsi}.

\begin{lemma} \label{injectivityonboundary}
Let $E_0 > 0$ be given.
\begin{itemize}
\item[(i)] There exists a constant $\delta = \delta (E_0) > 0$ such that for all deformations
\begin{align}
\eta \in \mathcal{E} \quad \text{with} \quad \tilde{E}_{\operatorname{el}} (\eta) \leq E_0 \nonumber
\end{align}
it holds that
\begin{align}
X_1, X_2 \in \partial \Omega_0\quad \text{with} \quad \left| X_1 - X_2 \right| < \delta \quad \quad \Rightarrow \quad \quad \eta \left( X_1 \right) \neq \eta \left( X_2 \right). \nonumber
\end{align}
\item[(ii)] Let $\delta > 0$ be as in (i). Let in addition $\eta_0 \in \operatorname{int} (\mathcal{E})$ be given such that $\tilde{E}_{\operatorname{el}}(\eta_0) \leq E_0$ and $\eta_0$ is injective on $\partial \Omega_0$. Then there exists a constant $\Gamma = \Gamma (\eta_0,E_0)$ such that
\begin{align}
X_1, X_2 \in \partial \Omega_0 \quad \text{with} \quad \left| X_1 - X_2 \right| \geq \delta \quad \quad \Rightarrow \quad \quad \left| \eta \left( X_1 \right) - \eta \left( X_2 \right) \right| \geq \frac{\epsilon}{2} \nonumber
\end{align}

for all deformations
\begin{align}
\eta \in \mathcal{E} \quad \text{with} \quad \tilde{E}_{\operatorname{el}} (\eta) \leq E_0 \quad \text{and} \quad \left\| \eta - \eta_0 \right\|_{L^2(\Omega_0)} < \Gamma. \nonumber
\end{align}
\end{itemize}
\end{lemma}

\textbf{Proof}

For the proof of (i) we refer to \cite[Lemma 2.5]{variationalapproachtofsi}, a proof of the statement (ii) is given in \cite[Proposition 2.7]{variationalapproachtofsi}.
$\hfill \Box$\\

\section*{Acknowledgments}

This work has been supported by the Czech Science Foundation  through the Czech-Korean grant 
 22-08633J (\v{S}N) and the grant 23-04766S (BB \& JS) as well as the Ministry of Education, Youth and Sport of the Czech Republic through the grant LL2105 CONTACT (BB) and within the frame of project Ferroic Multifunctionalities (FerrMion) [project No. CZ.02.01.01/00/22\_008/0004591], within the Operational Programme Johannes Amos Comenius co-funded by
the European Union (JS). We also aknowledge the support of the Czech Academy of Sciences through the grant Praemium Academiae of {\v{S}}. Ne{\v{c}}asov{\'{a}} (\v{S}N \& JS). 
Finally, the Institute of Mathematics, CAS is supported by RVO:67985840.
 % The authors were supported by the Czech-Korean project GA\v{C}R/22-08633J.   
  
 %   {\v{S}}. Ne{\v{c}}asov{\'{a}}  and J. Scherz have been supported by the Czech-Korean project GA\v{C}R/ and  by  the Praemium Academiae of {\v{S}}. Ne{\v{c}}asov{\'{a}},


\begin{thebibliography}{9}

\bibitem[]{0}
\ 

\textcolor{blue}{\Large{Martin Kru\v{z}ik}}

\textcolor{red}{\textbf{Modeling}}

\bibitem{mk1}
E. Davoli, M. Kru\v{z}ík, P. Piovano and U. Stefanelli, \textit{Magnetoelastic thin films at large strains}, Continuum Mech. Thermodyn. 33: 327--341 (2021) 

essentially the same energy as in our case, but steady state case and in particular without any dissipation... also this is not an existence result, but it is the derivation of a model of a magnetoelastic film via dimension reduction ($\Gamma$-convergence) from a corresponding model in 3 dimensions

\bibitem{mk2}
S. Almi, M. Kru\v{z}ík and A. Molchanova, \textit{Linearization in magnetoelasticity}, arXiv:2401.09586

starting with the same nonlinear energy as for example in \cite{mk4} (with no regularization terms, admissible deformations are required to be injective and have strictly positive determinant of the gradient though), a linearized version of the model is derived via Gamma-convergence (the main result shows Gamma-convergence of the original energy to the linearized energy). There is no dissipation considered, the article only treats the steady-state case. In addition, the convergence of minimizers and minima of the nonlinear problem to minimizers and minima of the linearized problem is shown. The existence of such minimizers is not shown directly, however the existence of minimizers to the nonlinear (and thus the linearized) problem is already established for example in \cite{mk4}. \\

\textcolor{red}{\textbf{Existence}}

\bibitem{mk3}
M. Kru\v{z}ík, U. Stefanelli and J. Zeman, \textit{Existence results for incompressible magnetoelasticity}, Discrete Cont. Dyn. Systems 35: 2615--2623 (2015)

similar energy as in our case but without additional regularization terms, but slightly stronger assumptions on the elastic enery density, the material is incompressible and (most important) the dissipation is different from ours, it does not contain a transport term and apparently also not the time derivative of the deformation... the paper is mostly focused on the steady state case (without the dissipation) in which the existence of a minimizer to the energy is proved... for the full evolutionary system the existence of "energetic solutions" is then proved by a time discretization as in our case, where the existence of a minimizer in the discrete times follows from the steady state case

\bibitem{mk4}
M. Bresciani, E. Davoli and M. Kru\v{z}ík, \textit{Existence results in large-strain magnetoelasticity}, Ann. Inst. H. Poincaré Anal. Non Linéaire 40 no. 3: 557--592 (2023)

very similar to \cite{mk3}, but now for compressible instead of incompressible materials

\bibitem{mk5}
M. Kru\v{z}ík, U. Stefanelli and C. Zanini, \textit{Quasistatic evolution of magnetoelastic thin films via dimension reduction}, Disc. Cont. Dyn. Systems 35: 5999--6013 (2015)

quite similar to \cite{mk1}, again one starts with a three dimensional model and derives a model of a magnetoelastic film via dimension reduction. Here, however, \textcolor{red}{the elastic energy does not seem to be further speicified and I cannot see which regularizations are used?? Linear elasticity?}. One then proves the existence of "energetic solutions", again via time incremental minimization, of the bulk model and subsequently obtains a corresponding energetic solution to the model of the film by passing to the limit

\bibitem{mk6}
J. Ciambella, M. Kru\v{z}ík and G. Tomassetti, \textit{A theory of magneto-elastic nanorods obtained through rigorous dimension reduction}, Appl. Math. Modeling 106: 426--447 (2022)

Here a one dimensional rod model is derived via dimension reduction from a corresponding two dimensional problem. The energy is very close to our energy, including the second deformation gradient, the reciprocal of the determinant of the deformation gradient is also used, however not in the energy but in the set of admissible deformations. However, only the static case is considered and so there is no dissipation. The existence of a minimizer to the energy in the two-dimensional setting is shown via the direct method and, by a limit passage, a corresponding solution to the one-dimensional problem is obtained. \\

\textcolor{red}{\textbf{only elasticity, no micromagnetics} (I guess here we don't have to write too much because this problem is in particular entirely covered in \cite{variationalapproachtofsi})}

\bibitem{mk7}
\textit{On the passage from nonlinear to linearized viscoelasticity}

\bibitem{mk8}
\textit{Equilibrium of immersed hyperelastic solids}

\bibitem{mk13}
\textit{Equilibrium for multiphase solids with Eulerian interfaces}

\bibitem{mk9}
\textit{Separately global solutions to rate-independent processes in large-strain inelasticity}

\bibitem{mk10}
\textit{Linearization and Computation for Large-Strain Viscoelasticity}

\bibitem{mk11}
\textit{Nonlinear and linearized models in thermoviscoelasticity}

\bibitem{mk12}
\textit{Curvature-dependent Eulerian interfaces in elastic solids} \\

\textcolor{red}{\textbf{only micromagnetics, no elasticity}}

\bibitem{mk14}
M. Kru\v{z}ík and A.Prohl, \textit{Recent developments in the modeling, analysis, and numerics of ferromagnetism}, SIAM Review 48: 439--483 (2006)

different models are presented (the Landau-Lifshitz model in the stationary case and a simplified version thereof without the exchange energy as well as corresponding models (Landau-Lifshitz-Gilbert) in the dynamical setting) and existence results as well as numerical algorithms for these models are summarized... however, probably this is also not too interesting for us since it does not include any solid mechanics and the study of the Landau-Lifshitz-Gilbert equation in its own right is probably well-known

\bibitem{mk15}
\textit{Interactions between demagnetizing field and minor-loop development in bulk ferromagnets}

\bibitem{mk16}
\textit{Weierstrass-type maximum principle for microstructure in micromagnetics} \\

\textcolor{blue}{\Large{Tomá\v{s} Roubí\v{c}ek}}

\textcolor{red}{\textbf{Existence}}

\bibitem{tr1}
T. Roubí\v{c}ek, \textit{Landau theory for ferro-paramagnetic phase transition in finitely-strained viscoelastic magnets}, 2302.02850

A model for a combination of magnetoelasticity with thermodynamics is introduced and the existence of weak solutions to the model is proved. The model is formulated, as in Johannes' thesis, entirely in the current configuration. In particular, as opposed to our work, one does not solve for the deformation and the deformation does not change the shape of the material. The proof is also not obtained via time incremental minimization but via a Galerkin method. \textcolor{red}{However, the transport equation for the deformation gradient does not seem to be regularized here,} but a regularization of the energy via a second gradient ("hyper-stress") is used.

\bibitem{tr2}
T. Roubí\v{c}ek and G. Tomassetti, \textit{Phase transformations in electrically conductive ferromagnetic shape-memory alloys, their thermodynamics and analysis}, Archive Ration. Mech. Anal. 210: 1--43 (2013)

\textcolor{red}{A magneto-mechanical model is considered (the material is described as magnetostrictive but not explicitly magnetoelastic), additionally the material is electrically and thermally conducting. The magnetic force balance (see (30b)) seems to differ from our setting, apparently some form of the LLG equation); it does not contain a transport term. The existence of weak solutions to the model is proved via a time discretization and a regularization, the time discrete problem, however, is solved directly instead of via minimization}

\bibitem{tr3}
T. Roubí\v{c}ek and G.Tomassetti, \textit{A thermodynamically consistent model of magneto-elastic materials under diffusion at large strains and its analysis}, Zeit. angew. Math. Phys. 69 (2018)

The existence of weak solutions is proved for a model of magnetoelasticity in combination with thermodynamics and diffusion. The model ((3.19)-(3.21)) contains the same regularization in the equation of motion as our model (second gradient and a bound of the determinant away from zero) and in addition also takes inertia into account. The magnetic force balance instead is simpler than in our case (no transport term (just a simple time derivative) and no stray field). The proof does not make use of a time discretization at all, instead the result is proved via a Galerkin approximation and a regularization of the heat sources (Quote: "As we need to control the determinant of the deformation gradient, we cannot impose semi-convexity assumption and cannot rely on a time discretisation"). \\

\textcolor{red}{\textbf{only micromagnetics, no elasticity}}

\bibitem{tr4}
\textit{Microstructure evolution model in micromagnetics}

\bibitem{tr5}
\textit{Mesoscopic model for ferromagnets with isotropic hardening}

\bibitem{tr6}
\textit{A thermodynamically-consistent theory of the ferro/paramagnetic transition} \\

\textcolor{red}{\textbf{only elasticity, no micromagnetics} (I guess here we don't have to write too much because this problem is in particular entirely covered in \cite{variationalapproachtofsi})}

\bibitem{tr7}
T. Roubí\v{c}ek, \textit{Nonlinearly coupled thermo-visco-elasticity}, Nonlin. Diff. Eq. Appl. 20: 1243--1275 (2013)

Kind of similar to \cite{tr10}, here also the existence of weak solutions to a model for a combination of viscoelasticity with thermodynamics is proved. Again the elastic energy is regularized by the second deformation gradient (\textcolor{red}{I can't see boundedness of the determinant of the deformation gradient away from zero though}). Transport of the temperature is not considered, but inertia is considered as opposed to in \cite{tr10}. \textcolor{red}{Actually, I can't see see what the additional difficulties in \cite{tr10} compared to this article consist of.} The proof is also achieved via time-discretization and some regularization but the discretized system is solved directly instead of via minimization. 

\bibitem{tr8}
\textit{Rate-independent damage processes in nonlinear inelasticity}

\bibitem{tr9}
\textit{A note about hardening-free viscoelastic models in Maxwellian-type rheologies}

\bibitem{tr10}
A. Mielke and T. Roubí\v{c}ek, \textit{Thermoviscoelasticity in Kelvin-Voigt rheology at large strains},  Archive Ration. Mech. Anal. 238: 1--45 (2020)

\textcolor{red}{ARMA paper with A. Mielke} The existence of weak solutions to a model for a combination of viscoelasticity with thermodynamics is proved. The elastic part of the model uses a similar regularization as in our model (second deformation gradient and the reciprocal of the determinant of the deformation gradient). The dissipation potential does not seem to contain the transport of the temperature as it does not depend on the time derivative of the temperature. The proof is also achieved via time incremental minimization, but instead of one minimization problem one considers two minimization problems (one for the deformation and one for the temperature), which are solved iteratively (\textcolor{red}{Is this the "not full variational structure" of the proof?}). Moreover, a second approximation level is introduced, on which the system is regularized (\textcolor{red}{Is this the additional regularization needed due to the not full variational structure?})

\bibitem{tr11}
T. Roubí\v{c}ek, \textit{Thermodynamics of viscoelastic solids, its Eulerian formulation, and existence of weak solutions}, Zeitschrift f. angew. Math. Phys. 75 (2024)

The existence (and some regularity) of weak solutions to a model for a combination of viscoelasticity with thermodynamics is proved. Here, however, the model is formulated, as in Johannes' thesis and in [Landau theory for ferro-paramagnetic phase transition in finitely-strained viscoelastic magnets], entirely in the current configuration. In particular, as opposed to our work, one does not solve for the deformation and the deformation does not change the shape of the material. The proof is also not obtained via time incremental minimization but via a Galerkin method. \textcolor{red}{However, the transport equation for the deformation gradient does not seem to be regularized here,} but a regularization of the energy via a second gradient ("hyper-stress") is used.

\bibitem{tr12}
\textit{Dynamics of charged elastic bodies under diffusion at large strains}

\bibitem{tr13}
\textit{Magnetic shape-memory alloys: thermomechanical modeling and analysis}

\bibitem{tr14}
T. Roubí\v{c}ek, \textit{Interaction of finitely-strained viscoelastic multipolar solids and fluids by an Eulerian approach},  J. Math. Fluid Mech. 25 (2023)

fluid-structure interaction problem with viscoelastic solids. The system is formulated entirely in the current configuration and regularized via a second gradient, existence and regularity of weak solutions is proved.

\bibitem{tr15}
\textit{Quasistatic viscoelasticity with self-contact at large strains}

\end{thebibliography}

\begin{thebibliography}{9}

\bibitem{emek}
B.E. Abali and F.A. Reich, \textit{Thermodynamically consistent derivation and computation of
electro-thermo-mechanical systems for solid bodies}, Comput. Methods Appl. Mech. Engrg. 319:567–-595 (2017) 


\bibitem{alphonse}
A. Alphonse, D. Caetano, A. Djurdjevac and C.M. Elliott, \textit{Function spaces, time derivatives and compactness for evolving
families of Banach spaces with applications to PDEs}, Journal of Differential Equations 353: 268--338 (2023)

%\bibitem{friedrichsinequality}
%C. Amrouche and N. Seloula, \textit{Lp-theory for vector potentials and Sobolev's inequalities for vector fields: Application to the Stokes equations with nonstandard boundary conditions}, Mathematical Models and Methods in Applied Sciences 23 no. 1: 37--92 (2013)

%\bibitem{abf}
%P. Angot, C.H. Bruneau, and P. Fabrie, \textit{A penalization method to take into account obstacles in incompressible viscous flows}, Numer. Math. 81: 497--520 (1999)

\bibitem{antman}
S.S. Antman and R.C. Rogers, \textit{Steady-state problems of nonlinear electro-magneto-thermoelasticity}, Arch. Ration. Mech. Anal. 95: 279--323 (1986)

\bibitem{ausanio}
G. Ausanio, A.C. Barone, C. Hison, V. Iannotti, G. Mannara and L. Lanotte, \textit{Magnetoelastic sensor application in civil buildings monitoring}, Sensors and Actuators A: Physical 123--124: 290--295 (2005)

%\bibitem{bks}
%BENE\v{S}OVÁ, B., KAMPSCHULTE, M., SCHWARZACHER, S.: A variational approach to hyperbolic evolutions and fluid-structure interactions.  (2020)

%\bibitem{bmtfsi}
%L. B\v{a}lilescu, J.S. Martín and T. Takahashi, \textit{Fluid-rigid structure interaction system with Coulomb's law}, SIAM Journal on Mathematical Analysis 49 no. 6: 4625--4657 (2017)

\bibitem{barchiesi}
M. Barchiesi, D. Henao and C. Mora-Corral, \textit{Local invertibility in Sobolev spaces with applications to nematic elastomers and magnetoelasticity}, Archive for Rational Mechanics and Analysis 224 no. 2: 743--816 (2017)

\bibitem{johannesarticle}
B. Bene\v{s}ová, J. Forster, C. Liu, A. Schlömerkemper \textit{Existence of weak solutions to an evolutionary model for magnetoelasticity}, SIAM Journal on Mathematical Analysis 50 no. 1: 1200--1236 (2018)

\bibitem{variationalapproachtofsi}
B. Benešová, M. Kampschulte and S. Schwarzacher, \textit{A variational approach to hyperbolic evolutions and fluid-structure interactions}, J. Eur. Math. Soc. 26 no. 12: 4615--4697 (2023)

%\bibitem{incompressiblecase}
%B. Benešová, Š. Nečasová, J. Scherz and A. Schlömerkemper, \textit{Fluid-rigid body interaction in an incompressible electrically conducting fluid}, SIAM Journal on Mathematical Analysis 55 no. 2: 929--965 (2023)

\bibitem{bhattacharya}
K. Bhattacharya, \textit{Microstructure of Martensite: Why it forms and how it gives rise to the shape-memory effect}, Oxford University Press, Oxford 2003

\bibitem{bhattacharyajames}
K. Bhattacharya and R. D. James, \textit{The Material is the Machine}, Science 307:53--54 (2005)


%\bibitem{bielskigambin}
%W. Bielski and B. Gambin, \textit{Relationship between existence of energy minimizers of incompressible and nearly incompressible magnetostrictive materials}, Reports on Mathematical Physics 66 no. 2: 147--157 (2010)

%\bibitem{bienkowski}
%A. Bie\'{n}kowski and R. Szewczyk, \textit{New possibility of utilizing amorphous ring cores as stress sensor}, Physica Status Solidi (a) 189 no. 3: 787--790 (2002)

\bibitem{bienkowski2}
A. Bie\'{n}kowski and R. Szewczyk, \textit{The possibility of utilizing the high permeability magnetic materials in construction of magnetoelastic stress and force sensors}, Sensors and Actuators A: Physical 113 no. 3: 270--276 (2004)

%\bibitem{blancducomet}
%X. Blanc and B. Ducomet, \textit{Weak and strong solutions of equations of compressible magnetohydrodynamics}, Handbook of Mathematical Analysis in Mechanics of Viscous Fluids: 2869--2925, Springer, Cham (2016)

%\bibitem{fsiintro}
%T. Bodnár, G.P. Galdi, Š. Nečasová, \textit{Fluid-Structure Interaction and Biomedical Applications}, Birkhäuser, Basel (2014)

%\bibitem{boulakiaguerrero}
%M. Boulakia and S. Guerrero, \textit{A regularity result for a solid-fluid system associated to the compressible Navier-Stokes equations}, Ann. Inst. H. Poincaré Anal. Non Linéaire 26 no. 3: 777--813 (2009)

%\bibitem{cottetmaitre}
%C. Bost, G.H. Cottet and E. Maitre, \textit{Convergence Analysis of a Penalization Method for the Three-Dimensional Motion of a Rigid Body in an Incompressible Viscous Fluid}, SIAM J. Numer. Anal. 48: 1313--1337 (2010)

\bibitem{actuatorsandsensors}
J.R. Brauer, \textit{Magnetic Actuators and Sensors}, John Wiley \& Sons, Inc., Hoboken, New Jersey (2006)

\bibitem{bresciani}
M. Bresciani, E. Davoli and M. Kružík, \textit{Existence results in large-strain magnetoelasticity}, Ann. Inst. H. Poincaré Anal. Non Linéaire 40 no. 3: 557--592 (2023)

\bibitem{brown2}
W.F. Brown, \textit{Magnetoelastic Interactions}, Springer-Verlag, Berlin Heidelberg (1966)

\bibitem{brown}
W.F. Brown, \textit{Micromagnetics}, Interscience Publishers, New York (1963)

%\bibitem{bungartzschäfer}
%H.-J. Bungartz, M. Schäfer: \textit{Fluid-Structure Interaction: Modelling, Simulation, Optimisation}, Springer-Verlag, Berlin Heidelberg (2006)

%\bibitem{cabannes}
%H. Cabannes: \textit{Theoretical Magnetohydrodynamics}, Academic Press, New York (1970)

\bibitem{carbou}
G. Carbou, M.A. Efendiev and P. Fabrie, \textit{Global weak solutions for the Landau-Lifschitz equation with magnetostriction}, Mathematical Methods in the Applied Sciences 34 no. 10: 1274--1288 (2011)

\bibitem{antonin1}
A. Češík, G. Gravina, M. Kampschulte, \textit{Inertial evolution of non-linear viscoelastic solids in the face of (self-)collision}, Calc. Var. 63: 55 (2024)

\bibitem{antonin2}
A. Češík, G. Gravina, M. Kampschulte, \textit{Inertial (self-)collisions of viscoelastic solids with Lipschitz boundaries}, arXiv:2312.00431

\bibitem{chipot}
M. Chipot, I. Shafrir, V. Valente, and G. Vergara Caffarelli, \textit{On a hyperbolic-parabolic system arising in magnetoelasticity}, J. Math. Anal. Appl. 352 no. 1: 120--131 (2009)

\bibitem{ciarletnecas}
P.G. Ciarlet, J. Nečas, \textit{Injectivity and self-contact in nonlinear elasticity}, Archive for Rational Mechanics and Analysis 97 no. 3: 171--188 (1987)

%\bibitem{CST}
%C. Conca, J. San Martin and M. Tucsnak, \textit{Existence of solutions for the equations modelling the motion of a rigid body in a viscous fluid}, Commun. Partial Differential Equations 25: 1019--1042 (2000)

%\bibitem{cc}
%M. Coquerelle and G.H. Cottet, \textit{A vortex level set method for the two-way coupling of an incompressible fluid with colliding rigid bodies}, J. Comput. Phys. 227: 9121--9137 (2008)

%\bibitem{dautraylions}%modellingreference
%R. Dautray and P.L. Lions, \textit{Mathematical Analysis and Numerical Methods for Science and Technology. Volume 1: Physical Origins and Classical Methods}, Springer-Verlag, Berlin (2000)

%\bibitem{davidson}
%P.A. Davidson, \textit{An Introduction to Magnetohydrodynamics}, Cambridge University Press, United Kingdom (2001)

\bibitem{francescojoshua}
F. De Anna, J. Kortum, A. Schlömerkemper, \textit{Struwe-like solutions for an evolutionary model of magnetoviscoelastic fluids}, Journal of Differential Equations 309: 455--507 (2022)

\bibitem{degiorgi}
E. De Giorgi, \textit{New problems on minimizing movements}, Boundary Value Problems for PDE and Applications: 81--98 (1993)

%\bibitem{fractionaltraceinequality}
%L. Del Pezzo and J.D. Rossi, \textit{Traces for fractional Sobolev spaces with variable exponents}, Adv. Oper. Theory 2 no. 4:  435--446 (2017)


\bibitem{desimonedolzmann}
A. DeSimone and G. Dolzmann, \textit{Existence of minimizers for a variational problem in two‐dimensional nonlinear magnetoelasticity}, Arch. Rational Mech. Anal. 144: 107--120 (1998)

\bibitem{desimonejames}
A. DeSimone and R.D. James, \textit{A constrained theory of magnetoelasticity}, J. Mech. Phys. Solids 50 no. 2: 283--320 (2002)

\bibitem{desimone}
A. DeSimone, R.V. Kohn, S. Müller and F. Otto, \textit{Recent analytical developments in micromagnetics}, The science of Hysteresis: Physical Modeling, Micromagnetics and Magnetization Dynamics. Vol. 2: 269--381, Elsevier, Amsterdam (2006)

\bibitem{podio}
A. DeSimone and P. Podio-Guidugli, \textit{On the Continuum Theory of Deformable Ferromagnetic Solids}, Arch. Rational Mech. Anal. 136:201--233 (1996) 

%\bibitem{desjardinsesteban}
%B. Desjardins and M.J. Esteban, \textit{Existence of weak solutions for the motion of rigid bodies in a viscous fluid}, Arch. Rational Mech. Anal. 146: 59--71 (1999)

%\bibitem{compressibledesjardinsesteban}
%B. Desjardins and M.J. Esteban, \textit{On Weak Solutions for Fluid-Rigid Structure Interaction: Compressible and Incompressible Models}, Comm. in Partial Differential Equations 25 no. 7--8: 1399--1413 (2000)

%\bibitem{dipernalions}
%R.J. DiPerna and P.L. Lions, \textit{Ordinary differential equations, transport theory and Sobolev spaces}, Invent. Math. 98: 511--547 (1989)


\bibitem{dorfmann}
L. Dormann and R.W. Ogden, \textit{Nonlinear Theory of Electroelastic and Magnetoelastic Interactions}, Springer, New York (2014)


%\bibitem{dormyiskakov}
%DORMY, E., ISKAKOV, A.: On magnetic boundary conditions for non-spectral dynamo simulations. \textit{Geophy. Astrophy. Fluid Dyn.} \textbf{99-6}, 481-492 (2005)

%\bibitem{dreherjungel}
%M. Dreher and A. Jüngel, \textit{Compact families of piecewise constant functions in Lp(0,T;B)}, Konstanzer Schriften in Mathematik 292 (2011)

%\bibitem{du2}
%H. Du, Y. Shao and G. Simonett, \textit{On a thermodynamically consistent model for magnetoviscoelastic fluids in 3D}, Journal of Evolution Equations 24 (2024)

\bibitem{du}
H. Du, Y. Shao and G. Simonett, \textit{Well-posedness for magnetoviscoelastic fluids in 3D}, Nonlinear Analysis: Real World Applications 69: 103759 (2023)

\bibitem{ellahiani}
I. Ellahiani, E.--H. Essoufi and M. Tilioua, \textit{Global existence of weak solutions to a three-dimensional fractional model in magneto-viscoelastic interactions}, Boundary Value Problems 2017: 122 (2017)

\bibitem{eriksen}
J. Eriksen, \textit{On Formulating and Assessing Continuum Theories
of Electromagnetic Fields in Elastic Materials}, J Elasticity  87:95–-108 (2007)

\bibitem{eringenmaugin}
A.C. Eringen and G.A. Maugin, \textit{Electrodynamics of Continua II}, Springer-Verlag, New York (1990)

%\bibitem{evans}
%L.C. Evans, \textit{Partial Differential Equations}, American Mathematical Society, Providence, Rhode Island (1998)

%\bibitem{feireislbuch2}
%FEIREISL, E.: Dynamics of Viscous Compressible Fluids. Oxford University Press, Oxford (2004)

%\bibitem{feireisldynamics}
%E. Feireisl, \textit{Dynamics of viscous compressible fluids,} Oxford University Press, Oxford (2004)

%\bibitem{feireislcompactness}
%E. Feireisl, \textit{On compactness of solutions to the compressible isentropic Navier-Stokes equations when the density is not square integrable,} Comment. Math. Univ. Carolin. 42 no. 1: 83--98 (2001)

%\bibitem{feireisl}
%E. Feireisl, \textit{On the motion of rigid bodies in a viscous compressible fluid}, Arch. Rational Mech. Anal. 167: 281--308 (2003)

%\bibitem{incompressiblefeireisl}
%E. Feireisl, \textit{On the motion of rigid bodies in a viscous incompressible fluid}, J. Evol. Equ. 3: 419--441 (2003)

%\bibitem{compressibleviscousfluids}
%FEIREISL, E., KARPER, T., POKORNÝ, M.: Mathematical theory of compressible viscous fluids: analysis and numerics. Birkhäuser, Basel (2016)

%\bibitem{singularlimits}
%E. Feireisl and A. Novotný, \textit{Singular Limits in Thermodynamics of Viscous Fluids}, Birkhäuser, Basel (2017)

%\bibitem{feireisl3}
%E. Feireisl, A. Novotný and H. Petzeltová, \textit{On the domain dependence of solutions to the compressible Navier–Stokes equations of a barotropic fluid,} Math. Meth. Appl. Sci. 25: 1045--1073 (2002)

%\bibitem{feireisl2}
%E. Feireisl, A. Novotný and H. Petzeltová, \textit{On the Existence of Globally Defined Weak Solutions to the Navier-Stokes Equations,} J. Math. Fluid Mech. 3: 359--392 (2001)

%\bibitem{feireisl2}
%FEIREISL, E., NOVOTNÝ, A., PETZELTOVÁ, H.: On the Existence of Globally Defined Weak Solutions to the Navier-Stokes Equations. \textit{J. math. fluid mech.} \textbf{3}, 359-392 (2001)

\bibitem{johannesthesis}
J. Forster, \textit{Variational approach to the modeling and analysis of magnetoelastic materials}, Ph.D. thesis, Universität Würzburg (2016), urn:nbn:de:bvb:20-opus-147226

%\bibitem{francu}
%J. Franců, \textit{Weakly continuous operators. Applications to differential equations}, Appl. Math. 39: 45--56 (1994)

%\bibitem{diamagnetic}%modellingreference
%E.P. Furlani, \textit{Magnetophoretic separation of blood cells at the microscale}, Journal of Physics D: Applied Physics 40 no. 5: 1313--1319 (2007)

\bibitem{galdi}
G.P. Galdi, \textit{An Introduction to the Mathematical Theory of the Navier-Stokes Equations: Steady-State Problems}, Springer-Verlag, New York (2011)

%\bibitem{G2}
%G. P. Galdi, \textit{On the motion of a rigid body in a viscous liquid: A mathematical analysis with applications}, Handbook of Mathematical Fluid Dynamics, Vol. I: 655--791, Elsevier Sci, Amsterdam (2002)

\bibitem{gar07}
C.J. García-Cervera, \textit{Numerical micromagnetics: a review}, Bol. Soc. Esp. Mat. Apl. SeMA 39: 103--135 (2007)

\bibitem{sourav}
H. Garcke, P. Knopf, S. Mitra and A. Schlömerkemper, \textit{Strong well-posedness, stability and optimal control theory for a mathematical model for magneto-viscoelastic fluids}, Calc. Var. 61: 179 (2022)

%\bibitem{GGH13}
%M. Geissert, K. G\"otze and M. Hieber, \textit{Lp-theory for strong solutions to fluid-rigid body interaction in Newtonian and generalized Newtonian fluids}, Trans. Amer. Math. Soc. 365 no. 3: 1393--1439 (2013)

%\bibitem{gbl}
%J.F. Gerbeau, C. Le Bris and T. Lelièvre, \textit{Mathematical Methods for the Magnetohydrodynamics of Liquid Metals}, Oxford University Press, Oxford (2006)

% \bibitem{GH}
%D. G\' erard-Varet, M. Hillairet.
% \newblock Regularity issues in 
% the problem of fluid structure 
% interaction.
% \newblock{\em Archive for 
% Rational Mechanical Analysis}
% (2010); {\bf 195}, 2, 375--407.

%\bibitem{GH2} G\' erard-Varet D., M. Hillairet.
%\newblock Existence of weak solutions up to collision for viscous fluid-solid systems with slip.
%\newblock{\em Comm. Pure Appl. Math.} (2014);
%\textbf{67}, no. 12, 2022--2075.

% \bibitem{GHC} 
%D. G\' erard-Varet,  M. Hillairet  and  C. Wang,
% \newblock The influence of 
% boundary conditions on the 
% contact problem in a 3D 
% Navier-Stokes flow,
% \newblock{\em J. Math. Pures  Appl.} (2015);  \textbf{103}, no. 1, 1--38.

%\bibitem{gigli}
%N. Gigli and S.J.N. Mosconi, \textit{A variational approach to the Navier–Stokes equations}, Bulletin des Sciences Mathématiques 136 no. 3: 256--276 (2012)

\bibitem{gilbarg}
D. Gilbarg and N.S. Trudinger, \textit{Elliptic Partial Differential Equations of Second Order}, Springer-Verlag, Berlin (2001)

%\bibitem{drugdelivery}
%T.S. Gregory, R. Cheng, G. Tang, L. Mao and Z.T.H. Tse, \textit{The magnetohydrodynamic effect and its associated material designs for biomedical applications: a state-of-the-art review}, Adv. Funct. Mater. 26: 3942--3952 (2016)

%\bibitem{endocapsules}
%T.S. Gregory, K.J. Wu, J. Yu, J.B. Box, R. Cheng, L. Mao, G. Tang and Z.T.H. Tse, \textit{Magnetohydrodynamic-Driven Design of Microscopic Endocapsules in MRI}, IEEE/ASME Transactions on Mechatronics 20 no. 6: 2691--2698 (2015)

%\bibitem{greinere} %modellingreference
%W. Greiner, \textit{Classical Electrodynamics}, Springer-Verlag, New York (1998)

%\bibitem{griffiths}%modellingreference
%D.J. Griffiths, \textit{Introduction to Electrodynamics}, 3rd ed., Prentice Hall, Upper Saddle River (1999)

\bibitem{grimes}
C.A. Grimes, S.C. Roy, S. Rani and Q. Cai, \textit{Theory, instrumentation and applications of magnetoelastic resonance sensors: a
review}, Sensors 11 no.3: 2809--2844 (2011)

%\bibitem{guermondminev2d}
%J.L. Guermond and P.D. Minev, \textit{Mixed Finite Element Approximation of an MHD Problem Involving Conducting and Insulating Regions: The 2D Case}, ESAIM: Mathematical Modelling and Numerical Analysis 136 no. 33: 517--536 (2002)

%\bibitem{guermondminev}
%J.L. Guermond and P.D. Minev, \textit{Mixed Finite Element Approximation of an MHD Problem Involving Conducting and Insulating Regions: The 3D Case}, Numer. Meth. PDE 19: 709--731 (2003)

%\bibitem{GLSE}
%M.D. Gunzburger, H. C. Lee and A. Seregin, \textit{Global existence of weak solutions for viscous incompressible flow around a moving rigid body in three dimensions}, J. Math. Fluid Mech. 2: 219--266 (2000)

%\bibitem{haakmaitytakahashi}
%B. H. Haak, D. Maity, T. Takahashi and M. Tucsnak, \textit{Mathematical analysis of the motion of a rigid body in a compressible Navier-Stokes-Fourier fluid}, Math. Nachr. 292: 1972--2017 (2019)

%\bibitem{halphen}
%B. Halphen, Q.S. Nguyen, \textit{Sur les matériaux standard généralisés}, Journal de Mécanique 14 no. 1: 39--63 (1975)

\bibitem{healeykromer}
T.J. Healey and S. Krömer, \textit{Injective weak solutions in second-gradient nonlinear elasticity}, ESAIM: Control, Optimisation and Calculus of Variations 15 no. 4: 863--871 (2009)

\bibitem{stroffolini}
D. Henao and B. Stroffolini, \textit{Orlicz-Sobolev nematic elastomers}, Nonlinear Anal. 194: 111513 (2020)

%\bibitem{hl}
%M.I. Herreros and S. Ligüérzana, \textit{Rigid body motion in viscous flows using the finite element method}, Phys. Fluids 32 (2020)

%\bibitem{hiebermurata}
%M. Hieber and M. Murata, \textit{The Lp-approach to the fluid-rigid body interaction problem for compressible fluids}, Evol. Equ. Control Theory 4: 69--87 (2015)

%\bibitem{HOST}
%K.H. Hoffmann and V. N. Starovoitov, \textit{On a motion of a solid body in a viscous fluid. Two-dimensional case}, Adv. Math. Sci. Appl. 9: 633--648, (1999)

\bibitem{hs98}
A. Hubert and R. Schäfer, \textit{Magnetic Domains: The Analysis of Magnetic Microstructures}, Springer-Verlag, Berlin Heidelberg (1998)

%\bibitem{Imai}
%IMAI, I.: Chapter I. General principles of magneto-fluid dynamics. \textit{Suppl.
%Prog. Theor. Phys.} \textbf{24}, 1-34 (1962)

%\bibitem{jackson}
%J.D. Jackson, \textit{Classical Electrodynamics}, 3rd ed., John Wiley \& Sons, United States (1998)

\bibitem{kinderlehrer}
R.D. James and D. Kinderlehrer, \textit{Frustration and microstructure: an example in magnetostriction}, Prog. Part. Diff. Eqns.: calculus of variations, applications (Bandle, C. et. al., eds.) Pitman Res. Notes Math 267: 59--81 (1992)

\bibitem{jameskinderlehrer}
R.D. James and D. Kinderlehrer, \textit{Theory of magnetostriction with applications to $Tb_xDy_{1-x}Fe_2$}, Philos. Mag. 68: 237--274 (1993)

%\bibitem{jiliang}
%S. Jiang and F. Li, \textit{Convergence of the complete electromagnetic fluid system to the full compressible magnetohydrodynamic equations}, Asymptotic Analysis 95: 161--185 (2015)

%\bibitem{jiliang2}
%S. Jiang and F. Li, \textit{Rigorous derivation of the compressible magnetohydrodynamic equations from the electromagnetic fluid system}, Nonlinearity 25: 1735--1752 (2012)

%\bibitem{jiliang}
%S. Jiang and F. Li, \textit{Convergence of the complete electromagnetic fluid system to the full compressible magnetohydrodynamic equations}, Asymptotic Analysis 95: 161--185 (2015)

%\bibitem{jiliang2}
%S. Jiang and F. Li, \textit{Rigorous derivation of the compressible magnetohydrodynamic equations from the electromagnetic fluid system}, Nonlinearity 25: 1735--1752 (2012)

\bibitem{jiang}
N. Jiang, H. Liu and Y.-L. Luo, \textit{On well-posedness of an evolutionary model for magnetoelasticity: hydrodynamics of viscoelasticity and Landau-Lifshitz-Gilbert systems}, J. Differential Equations 367: 79--123 (2023)

\bibitem{joshuapaper}
M. Kalousek, J. Kortum and A. Schlömerkemper, \textit{Mathematical analysis of weak and strong solutions to an evolutionary model for magnetoviscoelasticity}, Discrete Contin. Dyn. Syst. Ser. S 14 no. 1: 17--39 (2021)

\bibitem{martinanja}
M. Kalousek and A. Schlömerkemper, \textit{Dissipative solutions to a system for the flow of magnetoviscoelastic materials}, Journal of Differential Equations 271: 1023--1057 (2021)

\bibitem{kovetz}
A. Kovetz, \textit{Electromagnetic Theory}, Oxford University Press (2000)
 
%\bibitem{joshuamastersthesis}
%J. Kortum, \textit{Existence of a solution for an evolutionary magnetostrictive model}, Master thesis, Universität Würzburg (2017)

\bibitem{joshuathesis}
J. Kortum, \textit{Global existence and uniqueness results for nematic liquid crystal and magnetoviscoelastic flows}, Ph.D. thesis, Universität Würzburg (2021), urn:nbn:de:bvb:20-opus-278271

%\bibitem{weakstronguniqueness}
%O. Kreml, \v{S}. Ne\v{c}asová and T. Piasecki, \textit{Weak-strong uniqueness for the compressible fluid-rigid body interaction}, J. Differential Equations 268 no. 8: 4756--4785 (2020)

\bibitem{kromerroubicek}
S. Krömer, T. Roubi\v{c}ek, \textit{Quasistatic viscoelasticity with self-contact at large strains}, Journal of Elasticity 142: 433--445 (2020)

\bibitem{kruzikprohl}
M. Kru\v{z}ík and A. Prohl, \textit{Recent developments in the modeling, analysis, and numerics of ferromagnetism}, SIAM Rev 48 no. 3: 439--483 (2006)

\bibitem{kruzikroubicek}
M. Kru\v{z}ík and T. Roubí\v{c}ek, \textit{Mathematical Methods in Continuum Mechanics of Solids}, Springer-Verlag, Cham (2019)

\bibitem{kruzikstefanellizeman}
M. Kru\v{z}ík, U. Stefanelli and J. Zeman, \textit{Existence results for incompressible magnetoelasticity}, Discrete and Continuous Dynamical Systems 35 no. 6: 2615--2623 (2015)

%\bibitem{kl}
%A.G. Kulikovskiy and G.A. Lyubimov, \textit{Magnetohydrodynamics}, Addison Wesley, Reading, Massachussets (1965)

%\bibitem{laforest}
%M. Laforest, \textit{The $p$-CurlCurl: Spaces, traces, coercivity and a Helmholtz decomposition in $L^p$.}

\bibitem{landaulifshitz}
L.D. Landau and E.M. Lifshitz, \textit{Electrodynamics of Continuous Media}, Pergamon Press, Oxford (1960)

\bibitem{elasticity}
L.D. Landau and E.M. Lifshitz, \textit{Theory of Elasticity}, Pergamon Press, Oxford (1970)

%\bibitem{optimalcontrol}
%LEE, E.B., MARKUS, L.: Foundations of Optimal Control Theory. John Wiley \& Sons, New York (1967)

%\bibitem{lions}
%P.L. Lions, \textit{Mathematical Topics in Fluid Mechanics: Volume 1: Incompressible Models}, Oxford University Press, Oxford (1996)

%\bibitem{compressiblelions}
%P.L. Lions, \textit{Mathematical Topics in Fluid Mechanics: Volume 2: Compressible Models}, Oxford University Press, Oxford (1998)

\bibitem{mielkeroubicek}
A. Mielke and T. Roubí\v{c}ek, \textit{Thermoviscoelasticity in Kelvin–Voigt Rheology at Large Strains}, Archive Ration. Mech. Anal. 238: 1–45 (2020)

%\bibitem{moffatt}
%MOFFATT, H.K.: Magnetic field generation in electrically conducting fluids, Cambridge University Press, Cambridge (1978)

%\bibitem{fsiintro2}
%Y. Modarres-Sadeghi, \textit{Introduction to Fluid-Structure Interactions}, Springer, Cham (2022)

%\bibitem{moreau}%modellingreference
%R. Moreau, \textit{Magnetohydrodynamics}, Springer, Dordrecht (1990)

%\bibitem{weakstronguniqueinc}
%B. Muha, \v{S}. Ne\v{c}asová and A. Rado\v{s}ević, \textit{A uniqueness result for 3D incompressible fluid-rigid body interaction problem}, J. Math. Fluid Mech. 23 no. 1: (2021)

%\bibitem{slipbcoffrictiontype}
%\v{S}. Ne\v{c}asová, J. Ogorzały, J. Scherz, \textit{The compressible Navier-Stokes equations with slip boundary conditions of friction type},  Z. Angew. Math. Phys. 74 no. 5: 188 (2023)

%\bibitem{compressiblecase}
%\v{S}. Ne\v{c}asová, M. Ramaswamy, A. Roy and A. Schlömerkemper, \textit{Motion of a rigid body in a compressible fluid with Navier-slip boundary condition}, J. Differential Equations 338: 256--320 (2022)

%\bibitem{alternativerigidbodies}
%\v{S}. Ne\v{c}asová, T. Takahashi and M. Tucsnak, \textit{Weak Solutions for the Motion of a Self-propelled Deformable Structure in a Viscous Incompressible Fluid}, Acta Appl. Math. 116: 329--352 (2011)

%\bibitem{nguyen}
%Q.S. Nguyen, \textit{Stability and Nonlinear Solid Mechanics}, Wiley, Chichester (2000)

\bibitem{novotnystraskraba}
A. Novotný and I. Stra\v{s}kraba, \textit{Introduction to the mathematical theory of compressible flow}, Oxford University Press, Oxford (2004)

%\bibitem{roberts}
%ROBERTS P.H.: An Introduction to Magnetohydrodynamics. Longmans, London (1967)

%\bibitem{robinson}
%ROBINSON J.C., RODRIGO J.L., SADOWSKI W.: The Three-Dimensional Navier-Stokes Equations: Classical Theory. Cambridge University Press (2016)

%\bibitem{rogers}
%ROGERS L.G.: Degree-independent Sobolev extension on locally uniform domains. \textit{Journal of Functional Analysis} \textbf{235}, 619-665 (2006)

\bibitem{ogden}
R.W.\ Ogden, \textit{Non-Linear Elastic Deformations}, Dover Publications Inc., New York (1997)

%\bibitem{oregan}
%D. O'Regan, \textit{Existence Theory for Nonlinear Ordinary Differential Equations}, Springer-Verlag, Dordrecht (1997)

\bibitem{frictionlesscontact1}
A.Z. Palmer and T.J. Healey, \textit{Injectivity and self-contact in second-gradient nonlinear elasticity}, Calc. Var. 56 no. 4: 114 (2017)

%\bibitem{dielectric}%modellingreference
%H. Pauly and H.P. Schwan, \textit{Dielectric properties and ion mobility in erythrocytes}, Biophys J. 6 no. 5: 621--639 (1966)

%\bibitem{milanpokornycompressiblense}
%M. Pokorný, \textit{Navier-Stokes equations}, Lecture notes, Charles University in Prague (2022), URL: \url{https://www.karlin.mff.cuni.cz/~pokorny/LectureNotes/NavierandStokes_eng.pdf} (Accessed: 22 September 2023)

\bibitem{praetorius}
D. Praetorius, \textit{Analysis of the operator $\Delta^{-1}$ div arising in magnetic models}, Zeitschrift für Analysis und ihre Anwendungen 23 no. 3: 589--605 (2004)

%\bibitem{rindler} F. Rindler, \textit{Calculus of Variations}, Springer, Cham (2018)

\bibitem{rogers}
R.C. Rogers, \textit{Nonlocal variational problems in nonlinear electromagneto-elastostatics}, SIAM Journal on Mathematical Analysis 19: 1329--1347 (1988)

\bibitem{tr1}
T. Roubí\v{c}ek, \textit{Landau theory for ferro-paramagnetic phase transition in finitely-strained viscoelastic magnets}, Mathematical Models and Methods in Applied Sciences 34 no. 2: 181-241 (2024)

\bibitem{roubicek}
T. Roubí\v{c}ek, \textit{Nonlinear Partial Differential Equations with Applications}, Birkhäuser Verlag, Basel, Boston, Berlin (2005) 

\bibitem{roubiceknew}
T. Roubí\v{c}ek, \textit{Visco-elastodynamics at large strains Eulerian}, Zeit. Angew. Math. Phys. 73: 80 (2022)

\bibitem{tr3}
T. Roubí\v{c}ek and G.Tomassetti, \textit{A thermodynamically consistent model of magneto-elastic materials under diffusion at large strains and its analysis}, Zeit. angew. Math. Phys. 69: 55 (2018)

\bibitem{roubicektomassetti}
T. Roubí\v{c}ek and G. Tomassetti, \textit{Phase transformations in electrically conductive ferromagnetic shape-memory alloys, their thermodynamics and analysis}, Archive Ration. Mech. Anal. 210: 1--43 (2013)

%\bibitem{roytakahashi}
%A. Roy and T. Takahashi, \textit{Stabilization of a rigid body moving in a compressible viscous fluid}, Journal of Evolution Equations 21: 167--200 (2020)

\bibitem{rybkaluskin}
P. Rybka and M. Luskin, \textit{Existence of energy minimizers for magnetostrictive materials}, SIAM Journal on Mathematical Analysis 36 no. 6: 2004--2019 (2005)

%\bibitem{bloodcompressible}
%A.H. Sacks and E.G. Tickner, \textit{The compressibility of blood}, Biorheology 5 no. 4: 271--274 (1968)

%\bibitem{tucsnak}
%J.A. San Martin, V. Starovoitov and M. Tucsnak, \textit{Global weak solutions for the two dimensional motion of several rigid bodies in an incompressible viscous fluid}, Arch. Rational Mech. Anal. 161: 93--112 (2002)

%\bibitem{sart}
%R. Sart, \textit{Existence of finite energy weak solutions for the equations MHD of compressible fluids}, Applicable Analysis 19 no. 3: 357--379 (2009)

%\bibitem{compressiblepaper}
%J. Scherz, \textit{Fluid-rigid body interaction in a compressible electrically conducting fluid}, Math. Nachr. (2023)

%\bibitem{thesis}
%J. Scherz, \textit{Weak Solutions to Mathematical Models of the Interaction between Fluids, Solids and Electromagnetic Fields}, PhD thesis, University of Würzburg and Charles University in Prague (2024)

\bibitem{thesis}
J. Scherz, \textit{Weak Solutions to Mathematical Models of the Interaction between Fluids, Solids and Electromagnetic Fields}, PhD thesis, University of Würzburg and Charles University in Prague, \url{urn:nbn:de:bvb:20-opus-349205} (2024)

\bibitem{zabensky}
A. Schlömerkemper, J. \v{Z}abenský, \textit{Uniqueness of solutions for a mathematical model for magneto-viscoelastic flows}, Nonlinearity 31 no. 6: 2989--3012 (2018)

%\bibitem{maxwellinequality}
%B. Schweizer: \textit{On Friedrichs inequality, Helmzholtz decomposition, vector potentials, and the div-curl lemma}, Trends on Applications of Mathematics to Mechanics - Springer INdAM series: 65--79, Springer, Cham (2018)

%\bibitem{stratton}
%STRATTON, J.A.: Electromagnetic Theory. McGraw-Hill Book Company, New York and London (1941)

%\bibitem{SST}
%San~Martin J. A., V. Starovoitov, M. Tucsnak.
%\newblock Global weak solutions for the two dimensional motion of several rigid bodies in an incompressible viscous fluid.
%\newblock {\em Arch. Rational Mech. Anal.},
%(2002); {\bf 161}, 93--112.

\bibitem{schloemerkemper}
A. Schl{\"o}merkemper, \textit{Mathematical Derivation of the
Continuum Limit of the Magnetic Force
between
Two Parts of a Rigid Crystalline}, Arch. Rational Mech. Anal. 176:227–-269 (2005) 

%\bibitem{bloodincompressible2}
%T.W. Secomb and A.R. Pries, \textit{Blood viscosity in microvessels: Experiment and theory}, Comptes Rendus Physique 14 no. 6: 470--478 (2013)

%\bibitem{SER3}
%D. Serre, \textit{Chute libre d'un solide dans un fluide visqueux incompressible. Existence}, Jap. J. Appl. Math. 4: 99--110 (1987)

\bibitem{micropumps}
A.R. Smith, A. Saren and J. Järvinen, K. Ullakko, \textit{Characterization of a high-resolution solid-state micropump that can be integrated into microfluidic systems}, Microfluid Nanofluid 18 no. 5: 1255--1263 (2015)

\bibitem{actuators}
R.L. Snyder, V.Q. Nguyen and R.V. Ramanujan, \textit{The energetics of magnetoelastic actuators is analogous to phase transformations in materials}, Acta Materialia 58 no. 17: 5620--5630 (2010)

%\bibitem{T} 
%T. Takahashi, \textit{Analysis of strong solutions for the equations modeling the motion of a rigid-fluid system in a bounded domain}, Adv. Differential Equations 8: 1499--1532 (2003)

%\bibitem{temam}
%R. Temam, \textit{Navier-Stokes Equations: Theory and Numerical Analysis}, North-Holland Publishing Company, Amsterdam (1979)

\bibitem{tiersten}
H.F. Tiersten, \textit{Coupled magnetomechanical equations for magnetically saturated insulators}, J. Math. Phys. 5: 1298--1318 (1964)

%\bibitem{Wa}
%C. Wang, \textit{Strong solutions for the fluid-solid systems in a 2-D domain}, Asymptot. Anal. 89 no. 3--4: 263--306 (2014)

%\bibitem{webb}
%WEBB, J.R.L.: Extensions of Gronwall's inequality with quadratic growth terms and applications. \textit{Electronic Journal of Qualitative Theory of Differential Equations} \textbf{61}, 1-12 (2018)

\bibitem{xiaobhattacharya}
Y. Xiao and K. Bhattacharya, \textit{A Continuum Theory of Deformable, Semiconducting Ferroelectrics}, Arch. Rational Mech. Anal. 189:59--95 (2008)

\bibitem{yang}
Z. Yang and L. Zhang, \textit{Magnetic actuation systems for miniature robots: a review}, Advanced Intelligent Systems 2 no. 13: 2000082 (2020)

%\bibitem{zeidler}
%E. Zeidler, \textit{Nonlinear Functional Analysis and its Applications I: Fixed-Point Theorems}, Springer Verlag, New York (1986)

%\bibitem{ziemer}
%ZIEMER, W.P.: Weakly Differentiable Functions. Springer-Verlag, New York (1989)

\bibitem{zhao1}
W. Zhao, \textit{Local well-posedness and blow-up criteria of magneto-viscoelastic flows}, Disc. Con. Dyn. Sys. 38 no. 9: 4637--4655 (2018)

\bibitem{zhao2}
W. Zhao, \textit{Weak-strong uniqueness of incompressible magneto-viscoelastic flows}, Comm. Pure. Appl. Anal. 19 no. 5: 2907--2917 (2020)

%\bibitem{heartmonitoring}
%Y. Zhou, X. Zhao, J. Xu, Y. Fang, G. Chen, Y. Song, S. Li and J. Chen, \textit{Giant magnetoelastic effect in soft systems for bioelectronics}, Nature Materials 20 no. 12: 1670--1676 (2021)

\end{thebibliography}
\end{document}